\documentclass[smallextended]{svjour3}       
\smartqed  
\usepackage{graphicx}
\usepackage{amsfonts,anysize,amsmath,indentfirst,geometry}
\usepackage{bm}
\usepackage{lineno}
\usepackage[colorlinks,citecolor=blue,linkcolor=blue]{hyperref}
\usepackage{mathrsfs}
\usepackage{amssymb}
\usepackage{array}
\usepackage{booktabs}
\usepackage{multirow}
\usepackage{algorithm}
\usepackage{algorithmic}
\usepackage{hyperref}
\usepackage{float}
\usepackage{subfigure}
\usepackage[justification=centering]{caption}
\usepackage{marvosym}

\numberwithin{equation}{section}
\newtheorem{thm}{Theorem}
\newtheorem{lem}{Lemma}
\newtheorem{rmk}{Remark}
\newtheorem{prop}{Proposition}

\renewcommand{\algorithmicrequire}{\textbf{Input:}}
\renewcommand{\algorithmicensure}{\textbf{Output:}}

\DeclareGraphicsExtensions{.pdf,.jpeg,.jpg,.png,.eps}
\graphicspath{{./figures/},{./myfigures/}}

\begin{document}

\title{
Efficient Alternating Least Squares Algorithms for Low Multilinear Rank Approximation of Tensors}

\author{Chuanfu Xiao         \and
        Chao Yang \and
        Min Li 
}

\institute{Chuanfu Xiao \at
              School of Mathematical Sciences, Peking University, Beijing 100871, China \\
              \email{chuanfuxiao@pku.edu.cn}           
           \and
           Chao Yang (\Letter)\at
              School of Mathematical Sciences, Peking University, Beijing 100871, China\\
              \email{chao\_yang@pku.edu.cn}
           \and   
           Min Li \at
           Institute of Software, Chinese Academy of Sciences, Beijing 100190, China\\
           \email{limin2016@iscas.ac.cn}
}

\date{Received: date / Accepted: date}

\maketitle

\begin{abstract}
The low multilinear rank approximation, also known as the truncated Tucker decomposition, has been extensively utilized in many applications that involve higher-order tensors. Popular methods for low multilinear rank approximation usually rely directly on matrix SVD, therefore often suffer from the notorious intermediate data explosion issue and are not easy to parallelize, especially when the input tensor is large.
In this paper, we propose a new class of truncated HOSVD algorithms based on alternating least squares (ALS) for efficiently computing the low multilinear rank approximation of tensors. The proposed ALS-based approaches are able to eliminate the redundant computations of the singular vectors of intermediate matrices and are therefore free of data explosion. Also, the new methods are more flexible with adjustable convergence tolerance and are intrinsically parallelizable on high-performance computers. Theoretical analysis reveals that the ALS iteration in the proposed algorithms is q-linear convergent with a relatively wide convergence region. Numerical experiments with large-scale tensors from both synthetic and real-world applications demonstrate that ALS-based methods can substantially reduce the total cost of the original ones and are highly scalable for parallel computing.
	
	\keywords{Low multilinear rank approximation \and Truncated Tucker decomposition\and Alternating least squares\and Parallelization}
	\subclass{15A69 \and 49M27 \and 65D15 \and 65F55}
\end{abstract}

\section{Introduction}
\label{S1}

As a natural extension of vectors (first-order) and matrices (second-order), higher-order tensors
have been receiving increasingly more attention in various applications, 
such as signal processing \cite{Lathauwer2004,Sidiropoulos2017}, 
computer vision \cite{Shashua2001,Vasilescu2002,Vlasi2005}, 
chemometrics \cite{Henrion1993,Jiang2000,Bro2006}, 
deep learning \cite{Charalampous2014,Novikov2015},
and scientific computing \cite{Beylkin2002,Khoromskij2012,Hackbusch2014}.
For decades, tensor decompositions have been extensively utilized
as an efficient tool for dimension reductions, latent variable analysis 
and other purposes in a wide range of scientific and engineering fields
\cite{Hackbusch2014,Beckmann2005,Kroonenberg2008,Ishteva2008,Holtz2012,Ding2016}.
There exist a number of tensor decomposition models, such as 
canonical polyadic (CP, or CANDECOMP/PARAFAC) decomposition \cite{Hitchcock1927,Carroll1970,Kiers2000},
Tucker decomposition \cite{Tucker1966,Levin1963,Lathauwer2000-1},
tensor train (TT) model \cite{Oseledets2011}, 
and hierarchical Tucker (HT) model \cite{Grasedyck2010,Hackbusch2009,Oseledets2009}.
Among them, the Tucker decomposition, also known as
the higher-order singular value decomposition (HOSVD),
is regarded as a generalization of the matrix singular value decomposition (SVD)
and has been applied with significant successes in many applications
\cite{Tucker1966,Levin1963,Lathauwer2000-1,Lathauwer2004,Ishteva2008,Kolda2009}.

In both theory and practice, a commonly considered tensor computation problem is the low multilinear rank approximation \cite{Lathauwer2000-2,Kolda2009}, also known as the truncated Tucker decomposition \cite{Vannieuwenhoven2012,Austin2016}, which reads
\begin{equation}\label{eq:prob}
\min_{\bm{\mathcal{B}}} \|\bm{\mathcal{A}} - \bm{\mathcal{B}}\|,
\end{equation}
where $\bm{\mathcal{A}}\in\mathbb{R}^{I_{1}\times I_{2}\times\cdots\times I_{N}}$ 
is a given $N$th-order tensor  and 
$\bm{\mathcal{B}}\in\mathbb{R}^{I_{1}\times I_{2}\times\cdots\times I_{N}}$
is its low multilinear rank approximation \cite{Lathauwer2000-1,Kolda2009}.
Existing approaches for solving (\ref{eq:prob}) can be roughly divided into 
two categories \cite{Vannieuwenhoven2011}: \emph{non-iterative} and \emph{iterative} methods.
The most popular \emph{non-iterative} algorithms for the low multilinear rank approximation
of higher-order tensors is the truncated HOSVD ($t$-HOSVD) \cite{Tucker1966,Lathauwer2000-2}
and its improved version, the sequentially truncated HOSVD ($st$-HOSVD) \cite{Vannieuwenhoven2012}.
Despite the fact that the results of $t$-HOSVD and $st$-HOSVD
are usually suboptimal, they can serve as good initial solution for popular iterative methods
such as higher-order orthogonal iteration (HOOI) \cite{Kroonenberg1980,Lathauwer2000-2}.
Other than the HOOI method, which is a first-order iterative method, some efforts are also made in developing
second-order approaches, such as Newton-type \cite{Elden2009,Savas2010,Ishteva2009}
and trust-region \cite{Ishteva2008,Ishteva2011} algorithms.
Although these methods can achieve faster convergence under certain conditions,
they are still in early study and are usually
not suitable for large-scale tensors \cite{Xu2018}.

In this paper, we focus on studying how to efficiently compute 
the truncated Tucker decomposition (\ref{eq:prob}) of higher-order tensors
by modifying the $t$-HOSVD and $st$-HOSVD algorithms so that certain fast iterative procedures are included.
As a major cost of the two algorithms, the computation of tensor-matrix multiplications
has been extensively optimized in a number of high-performance tensor libraries 
\cite{Rajbhandari2014,Li2015,Smith2015,Schatz2015}, and therefore is not the focus of this study.
Another potential bottleneck of the $t$-HOSVD and $st$-HOSVD algorithms
is the calculation of the singular vectors of the intermediate matrices, which can be done 
by applying SVD to the matricized tensor or eigen-decomposition to the Gram matrix
\cite{Kroonenberg1980,Lathauwer2000-2,Kolda2009,Austin2016,Oh2017}.
The matrix SVD can be obtained by using Krylov subspace methods 
\cite{Golub2008,Culum1983,Baglama2005}, whilst the eigen-decomposition of
the symmeric nonnegative definite Gram matrix can be done with a Krylov-Schur algorithm \cite{Sorensen1992,Stewart2001,Watkins2008,Golub2008}.
Due to the fact that these methods rely on the factorization of intermediate matrices,
they suffer from the notorious \emph{data explosion} issue \cite{Austin2016,Oh2017,Oh2019}.
And even if the hardware storage allowed, they are still not scalable for parallel computing
and the total computation cost could be unbearably high. 

In order to improve the performance of the $t$-HOSVD and $st$-HOSVD algorithms,
we propose a class of alternating least squares (ALS) based algorithms
for efficiently calculating the low multilinear rank approximation of tensors.
The key observation is that in the original algorithms 
the computations of singular vectors of the intermediate matrices
are indeed not necessary and can be replaced with low rank approximations,
and the low rank approximations can be done by using an ALS method which 
does not explicitly require intermediate tensor matricization, with the help of a row-wise update rule.
The proposed ALS-based algorithms enjoy advantages in computing efficiency,
error adaptivity and parallel scalability, especially for large-scale tensors.
We present theoretical analysis and show that the ALS iteration in the proposed algorithms
is q-linear convergent with a relatively wide convergence region. 
Several numerical experiments with both synthetic and real-world tensor data
demonstrate that new algorithms can effectively alleviate the data explosion issue
of the original ones and are highly parallelizable on parallel computers.

The organization of the paper is as follows. In Sec. \ref{S2}, we introduce some basic notations of tensor and 
the corresponding algorithms. In Sec. \ref{S3}, the $t$-HOSVD-ALS and $st$-HOSVD-ALS algorithms are proposed. 
Some theoretical analysis on the convergence behavior of the ALS methods can also be found in Sec. \ref{S3}.
After that, computational complexity and the approximation errors of proposed algorithms
are analyzed in Sec. \ref{S5}. Test results on several numerical experiments 
are reported in Sec. \ref{S6}. And the paper is concluded in Sec. \ref{S7}.

\section{Notations and Nomenclatures}\label{S2}
Symbols frequently used in this paper can be found in the following table.
\begin{table}[H]
	\normalsize
	\begin{center}
		\begin{tabular}{c|c}
			\hline
			Symbols   &  Notations \\
			\hline
			$a$ & scalar \\
			$\bm{a}$ & vector \\
			$\bm{A}$ & matrix  \\
			$\bm{\mathcal{A}}$ & three or higher-order tensor \\
			$\circ$ & vector outer product  \\
			$\times_{n}$ & mode-$n$ product of tensor and matrix \\
			$\bm{I}_{n}$ & identity matrix with size $n\times n$ \\
			$I_{n:N}$    & $\prod\limits_{i=n}^{N}I_{i}$ \\
			$\mathcal{R}(\bm{A})$ & a subspace formed by the columns of matrix $\bm{A}$ \\
			$\sigma(\bm{A})$ & a set that consists of sigular values of matrix $\bm{A}$ \\
			$\bm{A}^{\dagger}$ & pseudo-inverse of matrix $\bm{A}$\\
			\hline
		\end{tabular}
	\end{center}
\end{table}
Given an $N$th-order tensor 
$\bm{\mathcal{A}}\in\mathbb{R}^{I_{1}\times I_{2}\times\cdots\times I_{N}}$,
we denote $\bm{\mathcal{A}}_{i_{1},i_{2},\cdots,i_{N}}$
as its $(i_{1},i_{2},\cdots,i_{N})$-th element. In particular, rank one tensor is denoted as 
$$\bm{u}_{1}\circ\bm{u}_{2}\circ\cdots\circ\bm{u}_{N},$$
where $\bm{u}_{n}\in\mathbb{R}^{I_{n}}$ is a vector.

The Frobenius norm of tensor $\bm{\mathcal{A}}$ is defined as
$$\|\bm{\mathcal{A}}\|_{F} = \sqrt{\sum\limits_{i_{1},i_{2},\cdots,i_{N}}\bm{\mathcal{A}}_{i_{1},i_{2},\cdots,i_{N}}^{2}}.$$

The matricization of a higher-order tensor is a process
of reordering the elements of the tensor into a matrix.
For example, the mode-$n$ matricization of tensor $\bm{\mathcal{A}}$ is 
denoted as $\bm{A}_{(n)}$, which is a matrix belonging to
$\mathbb{R}^{I_{n}\times (I_{1}\cdots I_{n-1} I_{n+1}\cdots I_{N})}$.
Specifically, the $(i_{1},i_{2},\cdots,i_{N})$-th element of tensor $\bm{\mathcal{A}}$
is mapped to the $(i_{n},j)$-th entry of matrix $\bm{A}_{(n)}$, where
$$j = 1+\sum\limits_{k=1,k\neq n}^{N}(i_{k}-1)J_{k}\ \ \mbox{with}\ \ J_{k} = \prod\limits_{m=1,m\neq n}^{k-1}I_{m}.$$
The multilinear rank of a higher-order tensor $\bm{\mathcal{A}}$ is an integer array 
$(R_{1},R_{2},\cdots,R_{N})$, where $R_{n}$ is the rank of its
mode-$n$ matricization $\bm{A}_{(n)}$.

A frequently encountered operation in tensor computation
is the tensor-matrix multiplication. In particular,
the mode-$n$ tensor-matrix multiplication refers to the contraction of the tensor with a matrix
along the $n$-th index.
For example, suppose that $\bm{U}\in\mathbb{R}^{J\times I_{n}}$ is a matrix, 
the mode-$n$ product of $\bm{\mathcal{A}}$ and $\bm{U}$ is denoted as $\bm{\mathcal{A}}\times_{n}\bm{U}\in\mathbb{R}^{I_{1}\times\cdots\times I_{n-1}\times J\times I_{n+1}\times\cdots\ I_{N}}$. Elementwisely, one has
\begin{equation*}
\bm{\mathcal{B}}_{i_{1},\cdots,j,\cdots,i_{N}} = (\bm{\mathcal{A}}\times_{n}\bm{U})_{i_{1},\cdots,j,\cdots,i_{N}} = \sum\limits_{i_{n}=1}^{I_{n}}\bm{\mathcal{A}}_{i_{1},\cdots,i_{n},\cdots,i_{N}}\bm{U}_{j,i_{n}}.
\end{equation*}

The Tucker decomposition \cite{Tucker1966,Levin1963},
also known as the higher-order singular value decomposition (HOSVD) \cite{Lathauwer2000-1}, 
is formally defined as
\begin{equation*}
\bm{\mathcal{A}} = \bm{\mathcal{G}}\times_{1}\bm{U}^{(1)}\times_{2}\bm{U}^{(2)}\cdots\times_{N}\bm{U}^{(N)},
\end{equation*}
where $\bm{\mathcal{G}}\in\mathbb{R}^{R_{1}\times R_{2}\times\cdots\times R_{N}}$
is referred to as the core tensor
and $\bm{U}^{(n)}\in\mathbb{R}^{I_{n}\times R_{n}}$ are column orthogonal with each other 
for all $n\in\{1,2,\cdots,N\}$.
We remark here that the size of the core tensor is often smaller than that of 
the original tensor, though it is hard to know how small it can be a prior \cite{Silva2008,Kolda2009}.
In many applications, the Tucker decomposition is usually applied in its truncated form,
which reads 
\begin{equation}\label{Prob.2}
\begin{split}
\min\limits_{\bm{\mathcal{G}};\bm{U}^{(1)},\bm{U}^{(2)},\cdots,\bm{U}^{(N)}}\|\bm{\mathcal{A}} - \bm{\mathcal{G}}\times_{1}\bm{U}^{(1)}\times_{2}\bm{U}^{(2)}\cdots\times_{N}\bm{U}^{(N)}\|, \\
s.t.\ \ \ \ \bm{U}^{(n)T}\bm{U}^{(n)} = \bm{I}_{R_{n}},\ n\in\{1,2,\cdots,N\}
\end{split}
\end{equation}
where $(R_{1},R_{2},\cdots,R_{N})$ is a pre-determined truncation, smaller than 
the size of original tensor. 
Suppose that the exact solution of (\ref{Prob.2}) is
$\bm{U}^{*(1)}$, $\bm{U}^{*(2)}$, $\cdots$, $\bm{U}^{*(N)}$, and $\bm{\mathcal{G}}^{*}$,
then it is easy to see that
\begin{equation*}
\bm{\mathcal{G}}^{*} = \bm{\mathcal{A}}\times_{1}(\bm{U}^{*(1)})^{T}\times_{2}(\bm{U}^{*(2)})^{T}\cdots\times_{N}(\bm{U}^{*(N)})^{T},
\end{equation*}
which means 
\begin{equation}\label{eq:best}
\bm{\mathcal{A}}\times_{1}(\bm{U}^{*(1)})(\bm{U}^{*(1)})^{T}\times_{2}(\bm{U}^{*(2)})(\bm{U}^{*(2)})^{T}\cdots\times_{N}(\bm{U}^{*(N)})(\bm{U}^{*(N)})^{T}
\end{equation}
is the best low multilinear rank approximation of $\bm{\mathcal{A}}$. 

To compute the best low multilinear rank approximation of a higher-order tensor
in the truncated Tucker decomposition, a popular approach is 
the truncated HOSVD ($t$-HOSVD, \cite{Tucker1966}) originally presented by 
Tucker himself \cite{Tucker1966}. Nowadays, it is better known 
with the effort of Lathauwer \emph{et al.} \cite{Lathauwer2000-2}, who analyzed the structure of core tensor 
and proposed to employ SVD of the intermediate matrices in truncated HOSVD.
The computing procedure of  $t$-HOSVD is given in Algorithm \ref{A3.1}.
\begin{algorithm}[H] 
	\algsetup{linenosize=\normalsize}
	\normalsize 
	\caption{$t$-HOSVD \cite{Tucker1966,Lathauwer2000-2}}
	\label{A3.1}
	\begin{algorithmic}[1]
		\REQUIRE Tensor $\bm{\mathcal{A}}\in\mathbb{R}^{I_{1}\times I_{2}\times\cdots\times I_{N}}$, truncation $(R_{1},R_{2},\cdots, R_{N})$
		\ENSURE Low multilinear rank approximation $\hat{\bm{\mathcal{A}}}\approx\bm{\mathcal{G}}\times_{1}\bm{U}^{(1)}\times_{2}\bm{U}^{(2)}\cdots\times_{N} \bm{U}^{(N)}$
		\FORALL{$n\in\{1,2,\cdots,N\}$}
		\STATE $\bm{Q}\leftarrow$ leading left singular vectors of $\bm{A}_{(n)}$
		\STATE $\bm{U}^{(n)} \leftarrow \bm{Q}$
		\ENDFOR
		\STATE $\bm{\mathcal{G}} \leftarrow \bm{\mathcal{A}}\times_{1}\bm{U}^{(1)T}\times_{2}\bm{U}^{(2)T}\cdots\times_{N}\bm{U}^{(N)T}$
	\end{algorithmic}
\end{algorithm}

We remark here that in Algorithm \ref{A3.1}, the computation of $\bm{Q}$ can also be done by calculating
the $R_{n}$ eigenvectors of the Gram matrix $\bm{A}_{(n)}\bm{A}_{(n)}^{T}$.
It is clear that $t$-HOSVD can be seen as a natural extension of the truncated SVD of 
a matrix to higher-order tensors.
But unlike the matrix case, the approximation error of $t$-HOSVD 
is quasi-optimal \cite{Tucker1966,Lathauwer2000-2,Kolda2009,Vannieuwenhoven2012}.

As a subsequent improvement of $t$-HOSVD, the sequentially truncated HOSVD
($st$-HOSVD), proposed by Vannieuwenhoven \emph{et al.} \cite{Vannieuwenhoven2012}
uses a different truncation strategy, as shown in Algorithm \ref{A3.2}.
\begin{algorithm}[H]
	\algsetup{linenosize=\normalsize} 
	\normalsize
	\caption{$st$-HOSVD \cite{Vannieuwenhoven2012}}
	\label{A3.2}
	\begin{algorithmic}[1]
		\REQUIRE Tensor $\bm{\mathcal{A}}\in\mathbb{R}^{I_{1}\times I_{2}\times\cdots\times I_{N}}$, truncation $(R_{1},R_{2},\cdots, R_{N})$
		\ENSURE Low multilinear rank approximation $\hat{\bm{\mathcal{A}}}\approx\bm{\mathcal{G}}\times_{1}\bm{U}^{(1)}\times_{2}\bm{U}^{(2)}\cdots\times_{N} \bm{U}^{(N)}$
		\STATE Select an order of $\{1,2,\cdots,N\}$, i.e., $\{i_{1},i_{2},\cdots,i_{N}\}$. \\
		\STATE $\bm{\mathcal{B}}\ \leftarrow\ \bm{\mathcal{A}}$
		\FORALL{$n\in\{i_{1},i_{2},\cdots,i_{N}\}$}
		\STATE $\bm{U},\bm{\Sigma},\bm{V}^{T}\ \leftarrow$ matrix SVD of $\bm{B}_{(n)}$
		\STATE $\bm{U}^{(n)}\ \leftarrow\ \bm{U}$
		\STATE $\bm{\mathcal{B}}\ \leftarrow\ \bm{\Sigma}\bm{V}^{T}$ in tensor format
		\ENDFOR
		\STATE $\bm{\mathcal{G}}\ \leftarrow\  \bm{\mathcal{B}}$
	\end{algorithmic}
\end{algorithm}
Analogous to $t$-HOSVD, in Algorithm \ref{A3.2} one can also obtain $\bm{Q}$ by computing 
the $R_{n}$ eigenvectors of the Gram matrix $\bm{B}_{(n)}\bm{B}_{(n)}^{T}$, and the core tensor is updated by $\bm{\mathcal{B}} = \bm{\mathcal{B}}\times_{n}\bm{U}^{(n)T}$.  

Unlike $t$-HOSVD, the $st$-HOSVD algorithm does not compute the singular vectors of $\bm{A}_{(n)}$.
In particular, the factor matrices and core tensor in $st$-HOSVD are calculated simultaneously, which greatly reduces the size of the intermediate matrices.
It is worth mentioning that the order of $\{1,2,\cdots,N\}$ in Algorithm \ref{A3.2} could have a strong influence on 
the computational cost and approximation error of $st$-HOSVD. 
But there is no theoretical guidance on how to select the order. 
A heuristic suggestion on how to select the processing order of $st$-HOSVD was given based on the dimension of each mode \cite{Vannieuwenhoven2012}, i.e., $I_n$, $n=1,2,\cdots,N$.
As a supplement, the order can also be decided according to the truncation $(R_{1}$, $R_{2}$, $\cdots$, $R_{N})$ of the tensor, which is summarized in the following proposition.
\begin{prop}\label{prop1}
 Let $\bm{\mathcal{A}}\in\mathbb{R}^{I_{1}\times I_{2}\times\cdots\times I_{N}}$ be an $N$th-order tensor, and the truncation is set to $(R_{1}$, $R_{2}$, $\cdots$, $R_{N})$. 
 Without loss of generality, suppose that $I_{n}\approx I$ for any $n\in\{1,2,\cdots,N\}$, and $R_{1}\leq R_{2}\leq\cdots\leq R_{N}$. 
 Then the $st$-HOSVD algorithm based on order $\{1,2,\cdots,N\}$ has the lowest computational cost,
 as compared with other computational orders.
\end{prop}
\noindent{\bf\emph{Proof.}}
If we select $\{1,2,\cdots,N\}$ as the order of Algorithm \ref{A3.2}, when applying a Krylov subspace method to compute the truncated matrix SVD, the computational cost is 
\begin{equation}\label{Eq2.1}
\mathcal{O}(\sum\limits_{n=1}^{N}R_{1:n}I_{n:N}) \approx\mathcal{O}(\sum\limits_{n=1}^{N}R_{1:n}I^{N-n+1}),
\end{equation}
Similarly, the computational cost when we select $(i_{1},i_{2},\cdots,i_{N})$ as the order of Algorithm \ref{A3.2} is 
\begin{equation}\label{Eq2.2}
\mathcal{O}(\sum\limits_{n=1}^{N}R_{i_{1}:i_{n}}I_{i_{n}:i_{N}})\approx\mathcal{O}(\sum\limits_{n=1}^{N}R_{i_{1}:i_{n}}I^{N-n+1}).
\end{equation}
Clearly, (\ref{Eq2.1}) is smaller than (\ref{Eq2.2}).  \qed

For the best low multilinear rank approximation (\ref{Prob.2}),
it is easy to see is that $\bm{U}^{*(n)}$ is a column orthogonal factor matrix, 
therefore $(\bm{U}^{*(n)})(\bm{U}^{*(n)})^{T}$ represents the orthogonal projection 
of subspace $\mathcal{R}(\bm{U}^{*(n)})$. 
Consequently, subspace represented by the optimal factor matrices are critical. 
SVD and eigen-decomposition are the commonly applied approaches 
to determine this subspace in the original $t$-HOSVD and $st$-HOSVD procedures.
However, because of the introduction of the intermediate matrices,
both SVD and eigen-decomposition suffer from the notorious 
data explosion issue. Although some efforts have been made to alleviate 
the data explosion problem by, e.g., introducing an implicit Arnoldi procedure,
these fixes are usually not generalizable to large-scale tensors in real applications
and are not parallelization friendly.

In addition to SVD or eigen-decomposition, 
tensor matricization and tensor-matrix multiplication are also important 
in the original $t$-HOSVD and $st$-HOSVD algorithms. 
Recently, some efforts on high-performance optimizations of basic tensor operations 
are made. For example, Li \emph{et al.} proposed a shared-memory parallel implementation 
of dense tensor-matrix multiplication \cite{Li2015},  and
Smith \emph{et al.} considered sparse tensor-matrix multiplications \cite{Smith2015}.  
Nevertheless, the calculation of SVD or eigen-decomposition 
is still the major challenge in the $t$-HOSVD and $st$-HOSVD algorithms, especially for large-scale tensors.

\section{Alternating Least Squares Algorithms for $t$-HOSVD and $st$-HOSVD}\label{S3}

In this paper we tackle the challenges of the original
$t$-HOSVD and $st$-HOSVD algorithms from an alternating least squares (ALS) perspective.
Instead of utilizing SVD or eigen-decomposition on the intermediate matrices, we propose to compute
the dominant subspace with an ALS method to solve
a closely related matrix low rank approximation problem.
The classical ALS method for solving matrix low rank approximation problems
was originally proposed by Leeuw et al., \cite{Leeuw1976} and further applied in
principal component analysis \cite{Young1978}.
Algorithm \ref{A3.3} shows the detailed procedure of the ALS method.

\begin{algorithm}[H]
	\algsetup{linenosize=\normalsize}
	\normalsize
	\renewcommand{\algorithmicrequire}{\textbf{Input:}}
	\renewcommand{\algorithmicensure}{\textbf{Output:}}
	\caption{$[\bm{L}^{*},\bm{R}^{*}] = \mathrm{ALS}(\bm{A},r)$}
	\label{A3.3}
	\begin{algorithmic}[1]
		\REQUIRE Matrix $\bm{A}\in\mathbb{R}^{m\times n}$, truncation $r<\min\{m,n\}$ \\
		Initial guesses $\bm{L}_{0}\in\mathbb{R}^{m\times r}$ or $\ \bm{R}_{0}\in\mathbb{R}^{n\times r}$ \\
		\ENSURE  Low rank approximation $\hat{\bm{A}} = \bm{L}^{*} \bm{R}^{*T}$
		\STATE $k\ \leftarrow\ 0$
		\WHILE{not convergent}
		\STATE Solving multi-side least squares problem $\min\limits_{\bm{R}}\|\bm{L}_{k}\bm{R}^{T}-\bm{A}\|_{F}^{2}$
		\STATE $\bm{R}_{k}\ \leftarrow\  (\bm{A}^{T}\bm{L}_{k})(\bm{L}_{k}^{T}\bm{L}_{k})^{-1}$
		\STATE Solving multi-side least squares problem $\min\limits_{\bm{L}}\|\bm{R}_{k}\bm{L}^{T}-\bm{A}^{T}\|_{F}^{2}$
		\STATE $\bm{L}_{k+1}\ \leftarrow\ (\bm{AR}_{k})(\bm{R}_{k}^{T}\bm{R}_{k})^{-1}$
		\STATE $k\ \leftarrow\ k+1$
		\ENDWHILE
	\end{algorithmic}
\end{algorithm}

\begin{rmk}\label{rmk1}
We remark that unlike using an iterative method to solve matrix singular value problems, 
as was suggested in \cite{Vannieuwenhoven2012} to replace the matrix SVD,
the proposed ALS-based approach can avoid the computations of singular pairs/triplets and thus achieve higher performance.
\end{rmk}

\begin{rmk}\label{rmk2}
To consider the existence and uniqueness of the solution of the ALS iterations, we note that in line 3 of Algorithm \ref{A3.3}, if the coefficient matrices $\bm{L}_{k}$ are nonsingular, the multi-side least squares problem $\min\limits_{\bm{R}}\|\bm{L}_{k}\bm{R}^{T}-\bm{A}\|$ is equivalent to linear equation $\bm{L}_{k}^{T}\bm{L}_{k}\bm{R}^{T} = \bm{L}_{k}^{T}\bm{A}$, which has a unique solution as shown in line 4 of Algorithm \ref{A3.3}. The consistency of nonsingularity is further interpreted in Remark \ref{rmk3}.
\end{rmk}

As an iterative method, the number of iterations for the ALS method has a dependency on
the initial guess and the convergence criterion \cite{Szlam2017}.
In what follows we will establish a rigorous convergence theory of the ALS method
and derive an evaluation of the convergence region, which can help understand
how the initial guess could affect the speed of convergence.

To establish the convergence theory of the ALS method, we first require the following lemma,
which was proved in \cite{Wang2006}.
\begin{lem}\label{lem1}
	Let $\bm{A},\  \bm{B}\in\mathbb{R}^{n\times n}$ be symmetric positive definite matrices and satisfy
	$$\bm{B}\leq \bm{A},$$
	then the following inequalities hold
	$$\|\bm{A}^{-1}\bm{B}\|_{2}\leq1\,\,\mbox{and}\,\,\|\bm{BA}^{-1}\|_{2}\leq1,$$
	where $\bm{B}\leq\bm{A}$ represents $\bm{A} - \bm{B}$ is symmetric semi-positive matrix.
\end{lem}

The convergence theorem of Algorithm \ref{A3.3} is summarized in the theorem below.

\begin{thm}\label{thm2}	
	Let $\bm{A}\in\mathbb{R}^{m\times n}$ be a matrix, and $\sigma_{1}\geq\sigma_{2}\geq\cdots\geq\sigma_{\min\{m,n\}}$ be the singular values. Suppose that the following conditions hold:
	
	$1^{\circ}\ \sigma(\bm{L}_{k}),\ \sigma(\bm{R}_{k})$ are uniformly bounded. 
	
	$2^{\circ}\ \mathcal{R}(\bm{L}_{0})$ is in a neighborhood of the exact solution. 
	
	Then Algorithm \ref{A3.3} is local $q$-linear convergent, and the convergence ratio is approximately ${\sigma_{r+1}^{2}}/{\sigma_{r}^{2}}$, 
	where $\sigma_{r+1}<\sigma_{r}$.
\end{thm}

This theorem illustrates the convergence of the  ALS method in a viewpoint of subspace, 
and the convergence ratio depends on the gap of $\sigma_{r}$ and $\sigma_{r+1}$. 
The detailed proof can be found in Appendix~\ref{A}.

\begin{rmk}\label{rmk3}
	If condition $1^{\circ}$ in Theorem~\ref{thm2} is not satisfied, then either $\bm{L}_{k}$ or 
	$\bm{R}_{k}$ is close to singular. This implies that the truncation $r$ is inappropriately chosen, i.e., 
	greater than the numerical rank of $\bm{A}$. 
\end{rmk}

An evaluation of the convergence region of the ALS method can be found in the following theorem. 

\begin{thm}\label{thm3}
	Under the assumption of Theorem \ref{thm2}, provided that the initial guess $\bm{L}_{0}$ satisfies
	\begin{equation}\label{E2.8}
	\|\bm{L}_{0}^{(2)}(\bm{L}_{0}^{(1)})^{-1}\|_{2}\leq\sqrt{\frac{\sigma_{r}^{2}-(\sigma_{r}-\varepsilon)^{2}}{(\sigma_{r}-\varepsilon)^{2}-\sigma_{min}^{2}}},
	\end{equation}
	then the ALS method converges to the exact solution. Here 
	$$\bm{U}^{T}\bm{L}_{0} = \left(
	\begin{array}{c}
	\bm{U}_{1}^{T}\bm{L}_{0} \\
	\bm{U}_{2}^{T}\bm{L}_{0} \\
	\end{array}
	\right) = \left(
	\begin{array}{c}
	\bm{L}_{0}^{(1)} \\
	\bm{L}_{0}^{(2)} \\
	\end{array}
	\right),$$
	$\bm{A} = \bm{U}\bm{\Sigma}\bm{V}^{T}$ is the full SVD of $\bm{A}$, $\bm{U} = [\bm{U}_{1},\bm{U}_{2}]$ is the block form of $\bm{U}$, and $\varepsilon$ is an arbitary positive number such that
	$$\sigma_{r}-\varepsilon>\sigma_{r+1}.$$
\end{thm}

The proof of Theorem~\ref{thm3} is provided in Appendix~\ref{B}. 
We remark that it can be seen from the theorem that,
within the convergence region, a better initial guess guarantees faster convergence.
It is also worth noting that \eqref{E2.8} indicates that the convergence region depends on $\varepsilon$. 
A smaller $\varepsilon$ means higher requirement for the initial guess, but less number of iterations.

With the help of Algorithm \ref{A3.3}, we are able to solve the rank-$R_{n}$ approximation 
problem to obtain the dominant subspace of
$\bm{A}_{(n)}$ in $t$-HOSVD.
Based on it, we derive the ALS accelerated versions of the $t$-HOSVD algorithm,
namely $t$-HOSVD-ALS, presented in Algorithm \ref{A3.4}.

\vspace*{-1\baselineskip}
\begin{algorithm}[H]
	\algsetup{linenosize=\normalsize} 
	\normalsize
	\renewcommand{\algorithmicrequire}{\textbf{Input:}}
	\renewcommand{\algorithmicensure}{\textbf{Output:}}
	\caption{$t$-HOSVD-ALS}
	\label{A3.4}
	\begin{algorithmic}[1]
		\REQUIRE Tensor $\bm{\mathcal{A}}\in\mathbb{R}^{I_{1}\times I_{2}\times\cdots\times I_{N}}$, truncation $(R_{1},R_{2},\cdots, R_{N})$
		\ENSURE Low multilinear rank approximation $\hat{\bm{\mathcal{A}}}\approx\bm{\mathcal{G}}\times_{1}\bm{U}^{(1)}\times_{2}\bm{U}^{(2)}\cdots\times_{N} \bm{U}^{(N)}$
		\FORALL{$n\in\{1,2,\cdots,N\}$}
		\STATE $\bm{A}_{(n)}\leftarrow\bm{\mathcal{A}}$ in matrix format
		\STATE $\bm{L},\bm{R}\ \leftarrow\ \mathrm{ALS}(\bm{A}_{(n)},R_{n})$
		\STATE $ \hat{\bm{Q}},\hat{\bm{R}}\ \leftarrow$ reduced QR decomposition of $\bm{L}$
		\STATE $\bm{U}^{(n)}\ \leftarrow\ \hat{\bm{Q}}$
		\ENDFOR
		\STATE $\bm{\mathcal{G}}\ \leftarrow\ \bm{\mathcal{A}}\times_{1}\bm{U}^{(1)T}\times_{2}\bm{U}^{(2)T}\cdots\times_{N}\bm{U}^{(N)T}$
	\end{algorithmic}
\end{algorithm}
\vspace*{-0.5\baselineskip}

Similar to $st$-HOSVD, the proposed $t$-HOSVD-ALS algorithm does not explicitly compute the singular vectors of $\bm{A}_{(n)}$ either.
Instead, only an orthogonal basis is computed for the dominant subspace of $\bm{A}_{(n)}$ in $t$-HOSVD-ALS.
It is obtained by QR decomposition of $\bm{L}$, and the ALS method guarantees that $\mathcal{R}(\bm{L})$ is the left dominant subspace of $\bm{A}_{(n)}$.
In many applications the orthogonal basis suffices, but in case the singular vectors are required,
one can obtain them from the orthogonal basis by using, e.g., a low-overhead randomized approach \cite{Halko2011}.
Specifically, to calculate the singular vectors, line 5 of Algorithm \ref{A3.4} can be replaced with the following steps:
\begin{eqnarray*}
\bm{U},\bm{\Sigma},\bm{V}^{T} &\leftarrow& \mathrm{matrix\ SVD\ of\  }\hat{\bm{Q}}^{T}\bm{A}_{(n)},\\
\bm{U}^{(n)} &\leftarrow& \hat{\bm{Q}}\bm{U}.
\end{eqnarray*}

The ALS improved $st$-HOSVD algorithm, referred to as $st$-HOSVD-ALS,
can be analogously derived, as presented in Algorithm \ref{A3.5}. 

\vspace*{-1\baselineskip}
\begin{algorithm}[H]
	\algsetup{linenosize=\normalsize} 
	\normalsize
	\renewcommand{\algorithmicrequire}{\textbf{Input:}}
	\renewcommand{\algorithmicensure}{\textbf{Output:}}
	\caption{$st$-HOSVD-ALS}
	\label{A3.5}
	\begin{algorithmic}[1]
		\REQUIRE Tensor $\bm{\mathcal{A}}\in\mathbb{R}^{I_{1}\times I_{2}\times\cdots\times I_{N}}$, truncation $(R_{1},R_{2},\cdots, R_{N})$
		\ENSURE Low multilinear rank approximation $\hat{\bm{\mathcal{A}}}\approx\bm{\mathcal{G}}\times_{1}\bm{U}^{(1)}\times_{2}\bm{U}^{(2)}\cdots\times_{N} \bm{U}^{(N)}$
		\STATE Select an order of $\{1,2,\cdots,N\}$, i.e., $\{i_{1},i_{2},\cdots,i_{N}\}$
		\STATE  $\bm{\mathcal{B}}\ \leftarrow\ \bm{\mathcal{A}}$
		\FORALL{$n\in\{i_{1},i_{2},\cdots,i_{N}\}$}
		\STATE $\bm{B}_{(n)}\leftarrow\bm{\mathcal{B}}$ in matrix format 
		\STATE $\bm{L},\bm{R}\ \leftarrow\ \mathrm{ALS}(\bm{B}_{(n)},R_{n})$
		\STATE $ \hat{\bm{Q}},\hat{\bm{R}}\ \leftarrow$ reduced QR decomposition of $\bm{L}$
		\STATE $\bm{U}^{(n)}\ \leftarrow\ \hat{\bm{Q}}$
		\STATE $\bm{B}_{(n)}\ \leftarrow\ \hat{\bm{R}}\bm{R}^{T}$
		\STATE $\bm{\mathcal{B}}\leftarrow \bm{B}_{(n)}$ in tensor format
		\ENDFOR
		\STATE $\bm{G}_{(i_{N})}\ \leftarrow\ \bm{B}_{(i_{N})}$ 
		\STATE $\bm{\mathcal{G}}\leftarrow \bm{G}_{(i_{N})}$ in tensor format
	\end{algorithmic}
\end{algorithm}

The difference between Algorithm~\ref{A3.4} and \ref{A3.5} is whether or not
to store $\bm{R}$ and $\hat{\bm{R}}$, the right factor matrices of the ALS method
and reduced QR decomposition, respectively. 
Storing them will help reduce the overall computational cost when 
updating tensor $\bm{\mathcal{B}}$, and core tensor $\bm{\mathcal{G}}$ can be calculated with the last factor matrix simultaneously in Algorithm \ref{A3.5}.
Apart from the computational cost of ALS in the $t$-HOSVD-ALS algorithm, 
calculating the core tensor $\bm{\mathcal{G}}$ is also critical, especially for higher-order tensors. 

It is worth mentioning that the matricizations of $\bm{\mathcal{A}}$ and $\bm{\mathcal{B}}$ are not necessary if a row-wise update rule is used in $t$-HOSVD-ALS and $st$-HOSVD-ALS algorithms. Taking $t$-HOSVD-ALS as an example, in line 3 of Algorithm \ref{A3.4} we calculate the factor matrix $\bm{U}^{(n)}$ by using an ALS method to obtain the rank-$R_{n}$ approximation of $\bm{A}_{(n)}$. The key computation is solving a multi-side least squares problem with the right-hand side $\bm{A}_{(n)}$ or $\bm{A}_{(n)}^{T}$. And the multi-side least squares problem is equivalent to a series of independent least squares problems whose right-hand sides are the columns of $\bm{A}_{(n)}$ or $\bm{A}_{(n)}^{T}$, i.e., the mode-$n$ fiber or slice of tensor $\bm{\mathcal{A}}$. These independent least squares problems can be solved in parallel, and each least squares problem only requires a fiber or slice of tensor $\bm{\mathcal{A}}$, therefore the explicit matricization in line 2 of Algorithm \ref{A3.4} can be naturally avoided.

Compared with $t$-HOSVD and $st$-HOSVD, the proposed algorithms exhibit several advantages. 
First, the redundant computations of the singular vectors are totally avoided, 
thus the overall cost of the algorithm can be substantially reduced.
Second, the convergence of the ALS procedure is controllable by adjusting
the convergence tolerance. This is helpful considering the fact that
$t$-HOSVD and $st$-HOSVD are quasi-optimal, and are often used 
as the initial guess for other iterative algorithms such as HOOI.
Third, the algorithms are free of intermediate data explosion since the least square
problems can be solved without any intermediate matrices.

An added benefit of the proposed $t$-HOSVD-ALS and $st$-HOSVD-ALS algorithms is
that the solution of the multi-side least squares problems is intrinsically parallelizable.
By using the ALS method, each row of the factor matrix $\bm{L}$ or $\bm{R}$ can be independently updated.
Therefore, one can distribute the computation of the rows over multiple computing units. 
Since the workload for each row is almost identical, a simple static load distribution strategy suffices. 
All other operations in the algorithms, such as the matrix-matrix multiplication, 
the QR reduction and the small-scale matrix inversion, can also be easily parallelized by calling 
vendor-supplied highly optimized linear algebra libraries.  

\section{Computational Cost and Error Analysis}\label{S5}
In the proposed $t$-HOSVD-ALS and $st$-HOSVD-ALS algorithms, the performance of the ALS iteration depends on several factors, such as the initial guess and the convergence criterion.
Based on the convergence property and the convergence condition of the ALS method, 
we suggest to set the initial guess $\bm{L}_{0}$ as follows.
\begin{enumerate}
	\item Generate a random matrix $\bm{S}$, whose entries are uniform distributions on interval $[0,1]$.
	\item Compute the reduced QR decomposition $\bm{A}_{(n)}\bm{S}=\bm{QR}$.
	\item Let $\bm{Q}$ be the initial guess, i.e., $\bm{L}_{0} = \bm{Q}$.
\end{enumerate}
In this way, it is assured that $\mathcal{R}(\bm{L}_{0})$ is a subspace of $\mathcal{R}(\bm{A}_{(n)})$, 
which is closer to the left dominant subspace of $\bm{A}_{(n)}$ than a random initial guess. 
Also, step $3$ makes sure that the initial guess is properly normalized.
We remark here that unlike techniques utilized in randomized algorithms that require $\bm{S}$ to satisfy some special properties \cite{Halko2011,Mahoney2011,Mahoney2016},
the suggested initialization approach only requires $\bm{S}$ to be dense and full rank.

The stopping condition of the ALS iteration can be set to 
\begin{equation}\label{E5.0}
|\|\bm{A}_{(n)}-\bm{L}_{k}\bm{R}_{k}^{T}\|_{F}-\|\bm{A}_{(n)}-\bm{U}_{1}\bm{U}_{1}^{T}\bm{A}_{(n)}\|_{F}|\leq\eta {\|\bm{\mathcal{A}}\|_{F}},
\end{equation}
where $\mathcal{R}(\bm{U}_{1})$ is the left dominant subspace of $\bm{A}_{(n)}$, and $\eta$ is an accuracy tolerance parameter. 
In practice, however, $\bm{U}_{1}$ is often not available. 
We therefore advise to replace \eqref{E5.0} by 
\begin{equation}\label{E5.01}
|\|\bm{A}_{(n)}-\bm{L}_{k}\bm{R}_{k}^{T}\|_{F}-\|\bm{A}_{(n)}-\bm{L}_{k+1}\bm{R}_{k+1}^{T}\|_{F}|\leq\eta \|\bm{\mathcal{A}}\|_{F}
\end{equation}
as the stop criterion, which means that the relative approximation error almost no longer decreases, implying that $\{\bm{L}_{k},\bm{R}_{k}\}$ has converged to the critical point of optimization problem $\min\limits_{\bm{L},\bm{R}}\|\bm{A}_{(n)} - \bm{L}\bm{R}^{T}\|_{F}$.

Next, we will discuss truncation $R_{n}$ and how to select the tolerance parameter $\eta$ by error analysis.
To analyze the approximation error of ALS-based algorithms, we first recall a useful lemma.
\begin{lem}\label{lem2}\cite{Vannieuwenhoven2012}
	Let $\bm{U}^{(n)}\in\mathbb{R}^{I_{n}\times R_{n}}$, $n\in\{1,2,\cdots,N\}$ be a sequence of column orthogonal matrices, calculated via the $t$-HOSVD or $st$-HOSVD algorithm, and suppose that $\hat{\bm{\mathcal{A}}} = \bm{\mathcal{A}}\times_{1}(\bm{U}^{(1)}\bm{U}^{(1)T})\times_{2}(\bm{U}^{(2)}\bm{U}^{(2)T})\cdots\times_{N}(\bm{U}^{(N)}\bm{U}^{(N)T})$ is an approximation of $\bm{\mathcal{A}}\in\mathbb{R}^{I_{1}\times I_{2}\times\cdots\times I_{N}}$. Then 
	\begin{equation}\label{E5.1}
	\|\hat{\bm{\mathcal{A}}} - \bm{\mathcal{A}}\|_{F}^2
	\leq \sum\limits_{n=1}^{N}\gamma_{n}\leq
	N\|\bm{\mathcal{A}} - \bm{\mathcal{A}_{opt}}\|_{F}^2,
	\end{equation}
	where $\gamma_{n} = \sum\limits_{r = R_{n}+1}^{I_{n}}(\sigma_{r}^{(n)})^{2}$, and $\bm{\mathcal{A}_{opt}}$ is the optimal solution of problem (\ref{eq:prob}).
\end{lem} 

It is worth noting that although estimation \eqref{E5.1} ignores the computation error, it is still useful in practice. 
By Lemma \ref{lem2}, the error analysis of our algorithms is described in Theorem \ref{thm4}, with proof given in Appendix \ref{C}.

\begin{thm}\label{thm4}
	If the stop criterion of ALS is set to \eqref{E5.0}, then the approximation errors of Algorithm \ref{A3.4} and \ref{A3.5} are bounded by 
	\begin{equation}\label{E5.62}
	\frac{\|\hat{\bm{\mathcal{A}}} - \bm{\mathcal{A}}\|_{F}}{\|\bm{\mathcal{A}}\|_{F}}\leq \sqrt{\sum\limits_{n=1}^{N}(\eta_{n}^{2}+\frac{\gamma_{n}}{\|\bm{\mathcal{A}}\|_{F}^{2}})}\leq
	\sqrt{N}(\eta+\frac{\|\bm{\mathcal{A}} - \bm{\mathcal{A}_{opt}}\|_{F}}{\|\bm{\mathcal{A}}\|_{F}}),
	\end{equation}
	where $\eta = \max\limits_{n\in\{1,2,\cdots,N\}}\eta_{n}$.
\end{thm}

We remark that although  in practice \eqref{E5.0} is replaced by \eqref{E5.01},
numerical tests indicate that the main result \eqref{E5.62} still holds.
From \eqref{E5.62}, we advice to choose the tolerance parameter $\eta_{n}$ 
such that the dominant term in the right hand side of \eqref{E5.62} 
is ${\gamma_{n}}/{\|\bm{\mathcal{A}}\|_{F}^{2}}$ 
or ${\|\bm{\mathcal{A}} - \bm{\mathcal{A}_{opt}}\|_{F}}/{\|\bm{\mathcal{A}}\|_{F}}$.
Furthermore, if truncation $R_{n}$ is selected appropriately, both $\gamma_{n}$ and $\eta_{n}$
will be small, and the ALS will converge very fast since ${\sigma_{R_{n}+1}}/{\sigma_{R_{n}}}\ll1$.
On the other hand, less suitable truncation $R_{n}$ represents larger $\gamma_{n}$
and therefore larger $\eta_{n}$, which in turn reduces the required number of ALS iterations.  

Also of interest to us is the overall costs of the proposed algorithms.
We analyze cases related to both general higher-order tensor $\bm{\mathcal{A}}\in\mathbb{R}^{I_{1}\times I_{2}\times\cdots\times I_{N}}$ 
with truncation $(R_{1},R_{2},\cdots,R_{N})$ 
and cubic tensor $\bm{\mathcal{A}}\in\mathbb{R}^{I\times I\times\cdots\times I}$ 
with truncation $(R,R,\cdots,R)$.
The analysis results are shown in Table \ref{T5.1}, where $\mathrm{iter}_{n}$ is the number of ALS iterations for mode $n$. 
The proposed ALS-based algorithms, the computational costs rely greatly on $\mathrm{iter}_{n}$,
which depends on the initial guess, the truncation $R_{n}$ and the accuracy requirement. 
Our numerical result will reveal later that $\mathrm{iter}_{n}$ is usually far smaller than $R_{n}$,
which is consistent with previous studies of the ALS method for matrix computation \cite{Szlam2017}. 
\begin{table}[H]
	\begin{center}
		\caption{\normalsize Computational cost of different $t$- and $st$-HOSVD algorithms}
  \label{T5.1}
  \setlength{\tabcolsep}{5mm}{
  \begin{tabular}{c|c|c|c}
   \toprule
   \multicolumn{2}{c|}{Algorithms} &  $\bm{\mathcal{A}}\in\mathbb{R}^{I_{1}\times I_{2}\times\cdots\times I_{N}}$ &  $\bm{\mathcal{A}}\in\mathbb{R}^{I\times I\times\cdots\times I}$  \\
   \midrule
   \multicolumn{1}{c|}{\multirow{3}[8]*{$t$-HOSVD}} & ALS-based & $\mathcal{O}(\sum\limits_{n=1}^{N} (R_{n}I_{1:N})\mathrm{iter}_{n})$  & $\mathcal{O}(\sum\limits_{n=1}^{N} RI^{N}\mathrm{iter}_{n})$ \\
   \cline{2-4}
    & EIG-based & $\mathcal{O}(\sum\limits_{n=1}^{N}(I_{n}I_{1:N}+R_{n}I_{n}^{2}))$  & $\mathcal{O}(NI^{N+1}+NRI^{2})$ \\
    \cline{2-4}
    & SVD-based & $\mathcal{O}(\sum\limits_{n=1}^{N}R_{n}I_{1:N})$  & $\mathcal{O}(NRI^{N})$ \\
   \midrule
   \multicolumn{1}{c|}{\multirow{3}[8]*{$st$-HOSVD}} &ALS-based & $\mathcal{O}(\sum\limits_{n=1}^{N} (R_{1:n}I_{n:N})\mathrm{iter}_{n})$ & $\mathcal{O}(\sum\limits_{n=1}^{N}  (R^{n}I^{N-n+1})\mathrm{iter}_{n})$ \\
   \cline{2-4}
   &EIG-based & $\mathcal{O}(\sum\limits_{n=1}^{N} (R_{1:n-1}I_{n}I_{n:N}+R_{n}I_{n}^{2}))$ & $\mathcal{O}(\sum\limits_{n=1}^{N}R^{n-1}I^{N-n+2}+NRI^{2})$ \\
   \cline{2-4}
   &SVD-based & $\mathcal{O}(\sum\limits_{n=1}^{N}R_{1:n}I_{n:N})$ & $\mathcal{O}(\sum\limits_{n=1}^{N}R^{n}I^{N-n+1})$ \\
   \bottomrule
  \end{tabular}}
 \end{center}
\end{table}

For comparison purpose, we also list in the table the complexities of the original $t$- and $st$-HOSVD algorithms,
including the ones based on matrix SVD and those based on the eigen-decomposition of the Gram matrix.
From the table we can make the following observations.
\begin{itemize}
\item The EIG-based algorithms are usually most costly as compared to the SVD and ALS-based ones, 
especially when the truncation size is much smaller than the dimension length.
\item  It is hard to tell theoretically whether SVD-based algorithms are more  or less costly than the proposed the ALS-based algorithms.
We will examine their cost through numerical experiments in the next section.
\end{itemize}

\section{Numerical Experiments}\label{S6}
In this section, we will compare the proposed ALS-based algorithms with the original $t$-HOSVD and $st$-HOSVD algorithms by several numerical experiments related to both synthetic and real-world tensors. 
The implementation of the original $t$- and $st$-HOSVD algorithms includes \texttt{mlsvd} from Tensorlab \cite{Vervliet} and \texttt{hosvd} from Tensor Toolbox \cite{Bader}.
In particular, \texttt{mlsvd} utilizes matrix SVD and \texttt{hosvd} employs eigen-decomposition of Gram matrix, for computing the factor matrices, therefore we denote them as $t$($st$)-HOSVD-SVD and $t$($st$)-HOSVD-EIG, respectively.
To examine the numerical behaviors of these algorithms, we carry out most of the experiments in MATLAB R2019b on a computer equipped with an Intel Xeon Gold 6240 CPU of 2.60 GHz.
And to study the parallel performance of the proposed algorithms, we implement the algorithms 
in C++ and run them on a workstation equipped with an Intel Xeon Gold 6154 CPU of 3.00 GHz.
Unless mentioned otherwise, the tolerance parameter is set to $\eta=10^{-4}$, 
and the maximum number of ALS iterations is limited to $50$ in all tests.

\subsection{Reconstruction of random low-rank tensors with noise}\label{S6.1}
In the first set of experiments we examine the performance of the original $t$- and $st$-HOSVD algorithms,
as well as the proposed ALS-based ones for the reconstruction of random  low-rank tensors with noise.
The tests are designed following the work of \cite{Zhang2001,Vannieuwenhoven2011}.
Specifically, the input tensor is randomly generated as
\begin{equation*}
\hat{\bm{\mathcal{A}}} = \bm{\mathcal{A}} + \delta\,\bm{\mathcal{E}},
\end{equation*}
where the elements of $\bm{\mathcal{E}}$ follow the standard Gaussian distribution, 
and the noise level is controlled by $\delta = 10^{-4}$.
The base tensor $\bm{\mathcal{A}}$ is constructed by two models, i.e., CP model and Tucker model.
And to measure the reconstruction error  we use $\|\bm{\mathcal{B}} - \bm{\mathcal{A}}\|_{F}/\|\bm{\mathcal{A}}\|_{F}$, 
where $\bm{\mathcal{B}}$ is the low multilinear rank approximation of $\bm{\hat{\mathcal{A}}}$. 

\subsubsection{CP model}\label{S6.1.1}
First, we use the CP model to construct the base tensor $\bm{\mathcal{A}}\in\mathbb{R}^{I\times J\times K}$ as follows:
\begin{equation*}
\bm{\mathcal{A}} = \lambda_{1}\cdot \bm{a}_{1}\circ \bm{b}_{1}\circ \bm{c}_{1} + \lambda_{2}\cdot \bm{a}_{2}\circ \bm{b}_{2}\circ \bm{c}_{2} + \cdots + \lambda_{R}\cdot \bm{a}_{R}\circ \bm{b}_{R}\circ \bm{c}_{R},
\end{equation*}
where $\bm{a}_{r}\in\mathbb{R}^{I}$, $\bm{b}_{r}\in\mathbb{R}^{J}$, $\bm{c}_{r}\in\mathbb{R}^{K}$ are randomly generated normalized vectors,
and coefficients $\lambda_{r}\in[5,10]$ for all $r\in\{1,2,\cdots,R\}$.

\begin{figure}[!htb]
	\centering
	\includegraphics[width=0.99\hsize]{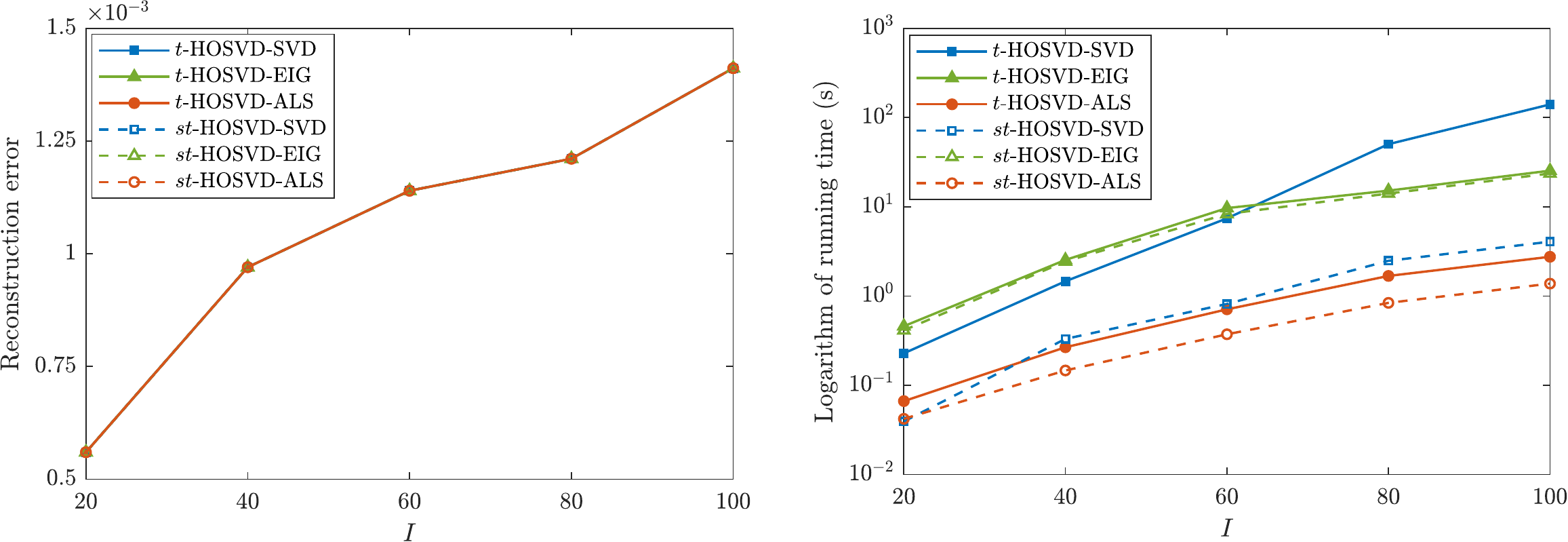}
	\caption{Reconstruction errors and running time of various low multilinear rank approximation algorithms for reconstructing random noisy low CP-rank tensors with gradually increased size.
	}\label{F6.1}
\end{figure}

In the experiments, we set the tensor size to be $J=I$ and $K=100I$ and
gradually increase $I$ from $20$ to $100$ with step $20$.
The truncation is set to $(R,R,R)$, where $R = 0.2 I$. 
We carry out the tests for 20 times and draw the averaged 
reconstruction errors and running time in Fig. \ref{F6.1}.
From the figure, it is observed that there is almost no difference in reconstruction error among all tested algorithms, indicating that the proposed ALS-based methods can maintain the accuracy of the original ones.
In terms of the running time, $t$-HOSVD-ALS is $3.4\times\sim 50.9\times$ faster than $t$-HOSVD-SVD, and $6.9\times\sim 13.6\times$ faster than $t$-HOSVD-EIG, respectively.
For $st$-HOSVD, the speedup of $st$-HOSVD-ALS is  $1.0\times\sim 3.0\times$ and 
$9.8\times\sim 22.3\times$, as compared to $st$-HOSVD-SVD and $st$-HOSVD-EIG, respectively.
We remark that the original algorithms behaves very differently in $t$-HOSVD and $st$-HOSVD.
In particular, changing from $t$-HOSVD-EIG to $st$-HOSVD-EIG leads to little performance improvement.
This is because the efficiency of eigen-decomposition of the Gram matrix
is strongly dependent on the size of the third mode of the input tensor,
and the sequentially updated algorithm is not able to help in this case.

\subsubsection{Tucker model}\label{S6.1.2}
Then we consider the base tesnor $\bm{\mathcal{A}}\in\mathbb{R}^{I_{1}\times I_{2}\times I_{3}\times I_{4}}$ constructed through the Tucker model in the following way:
$$\bm{\mathcal{A}} = \bm{\mathcal{G}}\times_{1}\bm{U}^{(1)}\times_{2}\bm{U}^{(2)}\times_{3}\bm{U}^{(3)}\times_{4}\bm{U}^{(4)},$$
where $\bm{\mathcal{G}}\in\mathbb{R}^{R_{1}\times R_{2}\times R_{3}\times R_{4}}$ is randomly generated core tensor whose elements follow the uniform distribution on interval $[5,10]$, and $\bm{U^{(i)}}\in\mathbb{R}^{I_{i}\times R_{i}},\ i=1,2,3,4$ are column orthogonal factor matrices.

\begin{figure}[!htb]
	\centering
	\includegraphics[width=0.99\hsize]{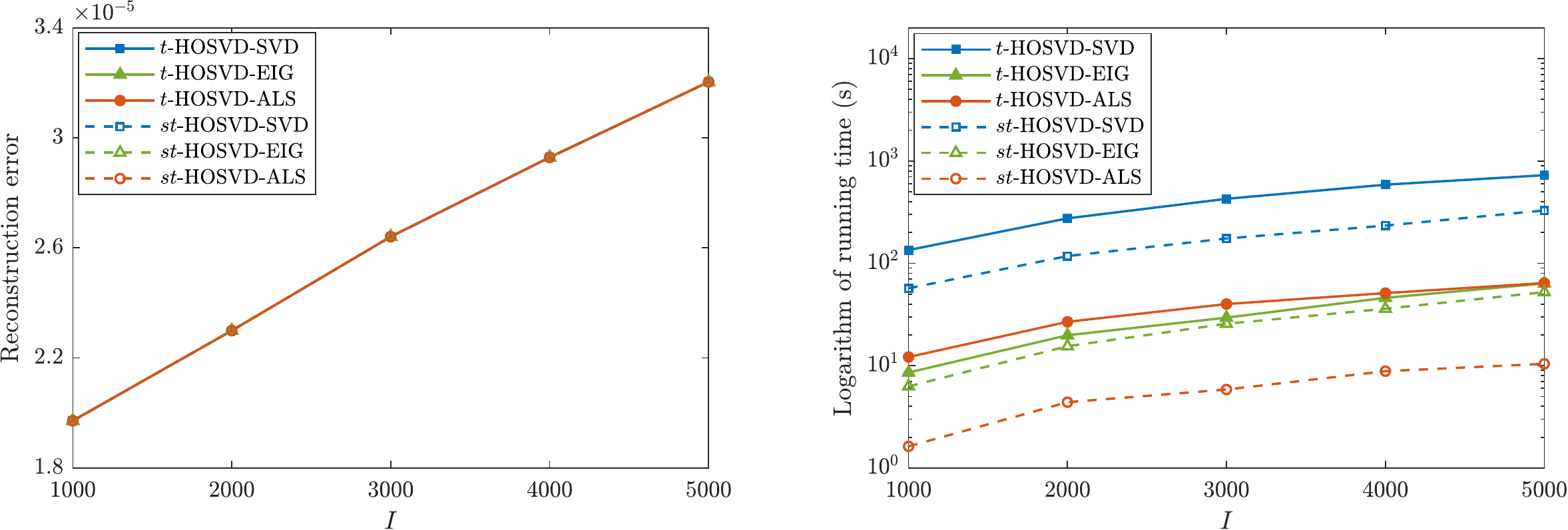}
	\caption{Reconstruction errors and running time of various low multilinear rank approximation algorithms for reconstructing random noisy low multilinear-rank tensors with gradually increased size.
	}\label{F6.1.2}
\end{figure}

In the experiments, we set $I_{2} = I_{3} = I_{4} = 100$ and gradually increase $I_{1} = I$ from $1,000$ to $5,000$ with step $1,000$. The truncation is set to be $(R,10,10,10)$, where $R = 0.01 I$. We again carry out the tests for 20 times, and draw the averaged reconstruction errors and running time in Fig.\ref{F6.1.2}. From the figure, we observe that the reconstruction errors of all tested algorithms are nearly identical. For $t$-HOSVD, the performance of $t$-HOSVD-ALS is close to that of $t$-HOSVD-EIG, and they are both $10.3\times$ faster than $t$-HOSVD-SVD. For $st$-HOSVD, the speedup of $st$-HOSVD-ALS is $26.5\times$ and $3.5\times$, as compared to $st$-HOSVD-SVD and $st$-HOSVD-EIG, respectively.

\subsection{Classification of handwritten digits}\label{S6.2}

The second set of experiments is designed for testing the capability of the original $t$- and $st$-HOSVD
algorithms and the proposed ALS-based ones on handwritten digits classification. 
It was studied that low multilinear rank approximation can be applied to compress the training data of images
so that the core tensor can be utilized for image classification on the test data \cite{Savas2007}.
In the tests, we use the MNIST database \cite{MNIST,LeCun1998} of handwritten digits.
We transfer the training dataset into a fourth-order tensor $\bm{\mathcal{A}}\in\mathbb{R}^{28\times28\times5000\times10}$, where the first- and second-mode are the texel modes, the third-mode corresponds to training images, and the fourth-mode represents image categories. In the tests, the truncation is fixed to be $(8,8,142,10)$, which is close to the setting of the reference work \cite{Savas2007}. We measure the classification accuracy of a classification algorithm as the percentage of the test images that is correctly classified.

\begin{table}[!htb]
 \normalsize
 \begin{center}
  \caption{\normalsize {Approximation error, classification accuracy and training time of various low multilinear rank approximation algorithms for handwritten digits classification of the MNIST database.}}\label{T6.2}
  \begin{tabular}{c|c|c|c|c|c|c}
   \toprule
   \multicolumn{1}{c}{\multirow{2}{*}{Algorithms}} & \multicolumn{3}{|c}{$t$-HOSVD} &\multicolumn{3}{|c}{$st$-HOSVD}\\
   \cline{2-7}
   & SVD & EIG & ALS & SVD & EIG &ALS \\
   \midrule
   Approximation error & $0.4330$ & $0.4330$ & $0.4335$ & $0.4288$ & $0.4288$ & $0.4291$\\
   \midrule
   Classification accuracy ($\%$)  & $95.45$ & $95.45$ & $95.45$ & $95.36$ & $95.36$ & $95.35$\\
   \midrule
   Training time (s) & $18.97$ & $5.38$ & $3.30$ & $4.07$ & $5.12$ & $1.53$ \\
   \bottomrule
  \end{tabular}
 \end{center}
\end{table}

We run the test for $20$ times and record the averaged approximation error, classification accuracy and running time of the tested algorithms in Table \ref{T6.2}. 
From the table it can be seen that the approximation errors and classification accuracies of all tested algorithms are almost indistinguishable, 
which again validates the accuracy of the proposed methods.
Table \ref{T6.2} also shows that the ALS-based methods are fastest in terms of training time. 
Specifically, the ALS-based approaches are on average $1.6\times\sim5.7\times$ and $2.7\times\sim3.3\times$ faster
than the original $t$-HOSVD and $st$-HOSVD algorithms, respectively.

\subsection{Compression of tensors arising from fluid dynamics simulations}\label{S6.3}
The purpose of this set of experiments is to examine the performance of different low multilinear rank approximation algorithms for compressing tensors generated from the simulation results of a lid-driven cavity flow, which is a standard benchmark for incompressible fluid dynamics \cite{Burggraf1966}.
The simulation is done in a square domain of length $1$ m with the speed of the top plate setting to $1$ m/s
and all other boundaries no flip. The kinematic viscosity is $\nu=1.0\times10^{-4}$ m$^2$/s, and the fluid properties is assumed to be laminar. We use the OpenFOAM software package \cite{Weller} to conduct the simulation
on a uniform grid with $100$ grid cells in each direction.
The simulation is run with time step $\Delta t = 1.0\times10^{-4}$ s and terminated at $t=1.0$ s. We record the magnitude of velocity at every time step. The simulation results of the lid-driven cavity flow are stored in a third-order tensor of size $100\times100\times10000$. To test the tensor approximation algorithms, we fix the truncation to $(20,20,20)$, corresponding to a compression ratio of $12,500:1$. 

\begin{figure}[!htb]
	\centering
	\includegraphics[width=0.55\hsize]{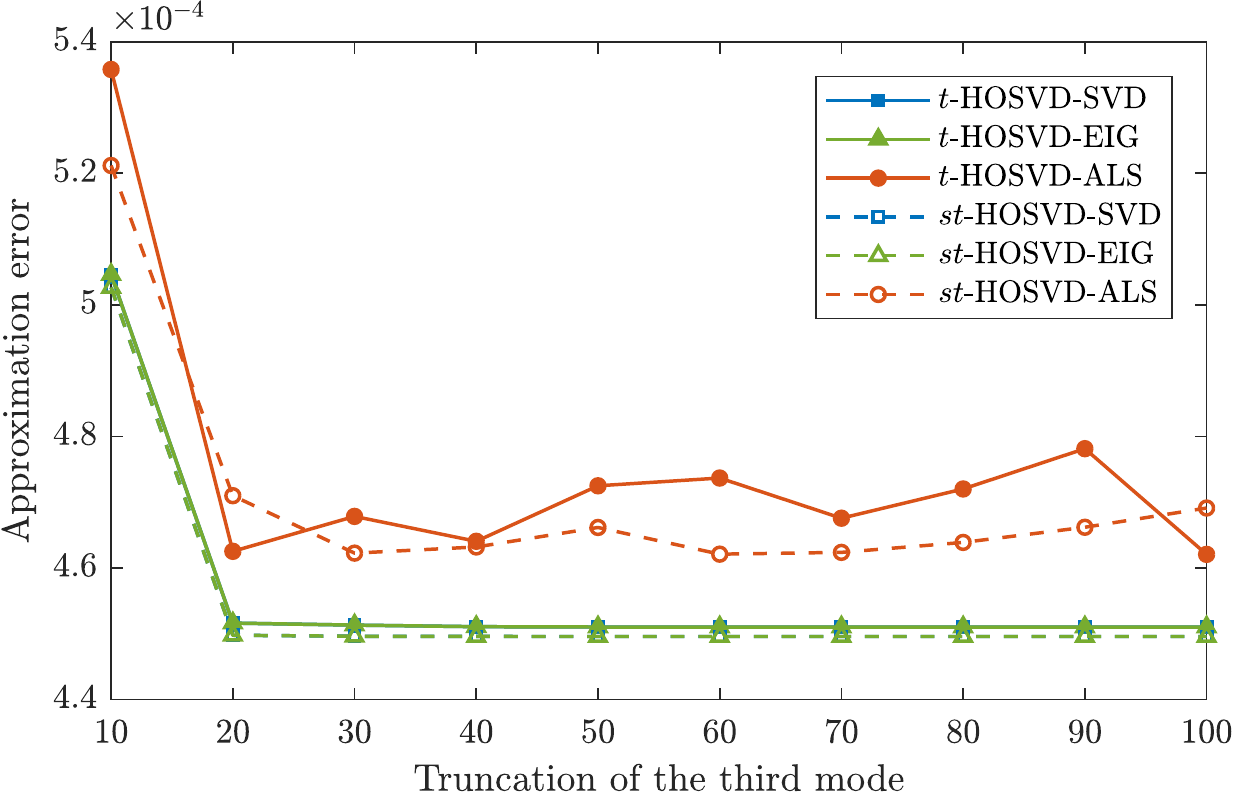}
	\caption{Approximation errors of $t$-HOSVD and $st$-HOSVD with gradually increased $R_3$.}\label{F6.3.2}
\end{figure}

First, we examine whether $(R_1, R_2, R_3) = (20,20,20)$ is a suitable truncation for the input tensor. Since the dimensions of the first two modes are both $100$, we consider $R_1=R_2=20$ is a relatively proper choice. And $R_3=20$ is remained to test for the third dimension of length $10,000$. Therefore, we perform a test to see how the approximation error varies as $R_3$ changes gradually from $10$ to $100$. The test results are shown in Fig. \ref{F6.3.2} , from which it can be observed that $R_3=20$ is also proper for the third dimension. 

\begin{table}[!]
	\normalsize
	\begin{center}
		\caption{\normalsize {Relative residuals and running time of various low multilinear rank approximation algorithms for compressing tensors arising from fluid dynamics simulations with different tolerance parameters.}}\label{T6.3-1}
		\begin{tabular}{c|c|c|c|c|c}
			\toprule
			\multicolumn{1}{c}{\multirow{2}{*}{$t$-HOSVD}}  &\multicolumn{1}{|c}{\multirow{2}{*}{SVD}} &\multicolumn{1}{|c}{\multirow{2}{*}{EIG}} & \multicolumn{3}{|c}{ALS} \\
			\cline{4-6}
			& & & $\eta = 10^{-2}$ & $\eta = 10^{-4}$ & $\eta = 10^{-6}$ \\
			\midrule
			Relative residual ($\times10^{-4}$) & $4.5161$ & $4.5161$& $15.7422$ & $4.6240$ & $4.5225$  \\
			\midrule
			Running time (s)  & $118.35$ & $23.61$ & $1.41+0.96$ & $2.21+0.96$& $3.964+0.96$ \\
			\bottomrule
			\multicolumn{1}{c}{\multirow{2}{*}{$st$-HOSVD}} &\multicolumn{1}{|c}{\multirow{2}{*}{SVD}} &\multicolumn{1}{|c}{\multirow{2}{*}{EIG}} & \multicolumn{3}{|c}{ALS} \\
			\cline{4-6}
			& & & $\eta = 10^{-2}$ & $\eta = 10^{-4}$ & $\eta = 10^{-6}$ \\
			\midrule
			Relative residual ($\times10^{-4}$) & $4.4976$ & $4.4976$& $14.2730$ & $4.7139$ & $4.5044$  \\
			\midrule
			Running time (s)  & $3.44$ & $21.70$ & $0.73$ & $1.16$& $1.82$ \\
			\bottomrule
		\end{tabular}
	\end{center}
\end{table}

Next, we study the efficiency of tested algorithms under different accuracy requirements. 
We run the test for 20 times with tolerance parameter $\eta$ adjusted to different values
and record the averaged relative residual and running time for each value of $\eta$; 
the test results are listed in Table~\ref{T6.3-1}.
Also provided in the table is the extra cost of computing the singular vectors for $t$-HOSVD-ALS, if requested.
From the table we have the following observations.
\begin{itemize}
	\item The relative residuals and running time of the original $t$- and $st$-HOSVD algorithms are independent of the change of the tolerance parameter $\eta$. This is due to the usage of Krylov subspace method for computing matrix truncated SVD or eigen-decomposition.
	\item For the ALS-based methods, the relative residuals and the running time both depend on $\eta$. 
	With $\eta$ decreased, the relative residuals are reduced to a similar level that the original algorithms can attain
	but more running time is required.
	\item The ALS-based algorithms are the fastest in all tests. It can achieve $6.1\times\sim45.6\times$ speedup for $t$-HOSVD
	and $2.3\times\sim18.4\times$ speedup for $st$-HOSVD, respectively.
	\item The overhead of computing the singular vectors for $t$-HOSVD-ALS is independent of the ALS tolerance and is relatively low.
\end{itemize}

\begin{figure}[!]
	\centering
	{\includegraphics[width=1\textwidth]{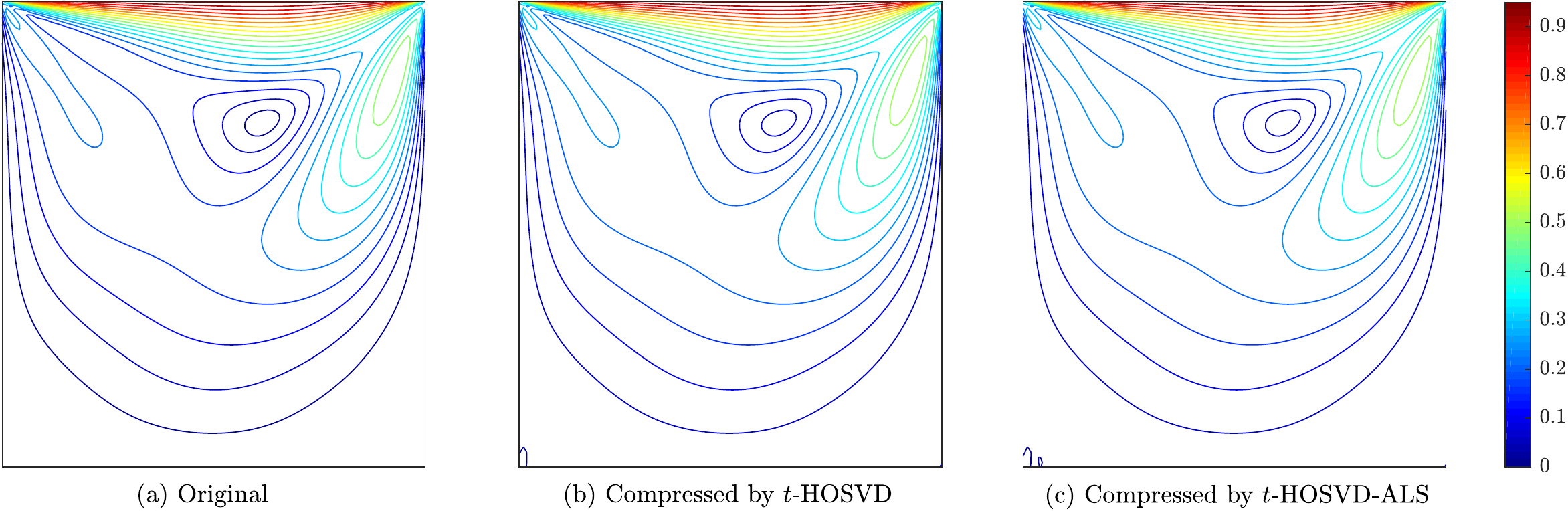}}
	\caption{\normalsize The original and compressed data at $t = 0.5$ s for compressing tensors arising from fluid dynamics simulations with tolerance parameter $\eta = 10^{-4}$.}\label{F6.3}
\end{figure}

It seems from the tests that despite the excellent performance of the proposed ALS-based methods,
the strong dependency between the sustained performance and the tolerance parameter $\eta$
could eventually lead to poor performance if $\eta$ is very small.
In practice the $t$- or $st$-HOSVD algorithm is often used as the initial guess of a supposedly more accurate
iterative method such as HOOI. 
In this case it is not necessary to use a very tight tolerance parameter.
To examine whether $\eta = 10^{-4}$ is a suitable choice for the ALS-based algorithms in the same tests,
we draw in Fig. \ref{F6.3} the contours of
the original and compressed velocity data at $t=0.5$ s. 
It clearly shows that when $\eta = 10^{-4}$, 
the compressed results are consistent with each other with very little discrepancy.
In fact, the measured maximum differences between the original data and the compressed data obtained by $t$-HOSVD and that by $t$-HOSVD-ALS 
are around $1.58\times10^{-3}$ and $1.33\times10^{-3}$, respectively,
which are very small considering that the compression ratio is over five orders of magnitude.

\begin{table}[!htb]
 \normalsize
 \begin{center}
 \caption{\normalsize {A comparison of HOOI results for compressing tensors arising from fluid dynamics simulations with initial solutions provided by various low multilinear rank approximation} algorithms.}\label{T6.3-2}
   \begin{tabular}{c|c|c|c}
   \toprule
   \multirow{2}{*}{Algorithm} & \multicolumn{2}{c|}{Relative residual} &\multirow{2}{*}{Number of HOOI iterations}\\
   \cline{2-3}
   &Initial  &Final & \\
   \midrule
   $t$-HOSVD-SVD & $4.5161\times10^{-4}$ &  $4.4807\times10^{-4}$  &  $9.0$  \\
   $t$-HOSVD-EIG & $4.5161\times10^{-4}$ &  $4.4807\times10^{-4}$  &  $7.0$  \\
   $t$-HOSVD-ALS & $4.6260\times10^{-4}$ &   $4.4807\times10^{-4}$  &  $6.0$  \\
   \midrule
   $st$-HOSVD-SVD & $4.4976\times10^{-4}$ & $4.4807\times10^{-4}$  &  $5.0$  \\
   $st$-HOSVD-EIG & $4.4976\times10^{-4}$ & $4.4807\times10^{-4}$  &  $6.0$  \\
   $st$-HOSVD-ALS & $4.5044\times10^{-4}$ &  $4.4807\times10^{-4}$  &  $7.0$  \\
   \bottomrule
  \end{tabular}
 \end{center}
\end{table}

To further investigate the applicability of the compressed results, we use the computed low multilinear rank approximation with tolerance parameter $\eta = 10^{-4}$ as the initial guess of the HOOI method with stopping criterion $10^{-12}$.
The HOOI method is obtained from the \texttt{tucker\_als} function of the Tensor Toolbox v3.1 \cite{Bader}. 
The test results are presented in Table \ref{T6.3-2}, 
in which we list the relative residuals with the $t$- and $st$-HOSVD provided initial guesses, the final relative residuals of HOOI, and the number of HOOI iterations, all averaged on 20 independent runs. From the table we can see that although the initial residual provided by the ALS-based algorithms 
are slightly larger than those provided by the original $t$- and $st$-HOSVD methods, same final residuals can be achieved after HOOI iterations nevertheless. And more importantly, the required numbers of HOOI iterations are insensitive to which specific $t$- and $st$-HOSVD algorithms, original or not, are used as shown in the tests. In other words, the proposed ALS-based methods are able to deliver similar results as the original ones when applying in HOOI, even when the tolerance parameter is relatively loose.

\subsection{Parallel performance}\label{S6.4}

\begin{figure}[!htb]
 \centering
 {\includegraphics[width=1.0\hsize]{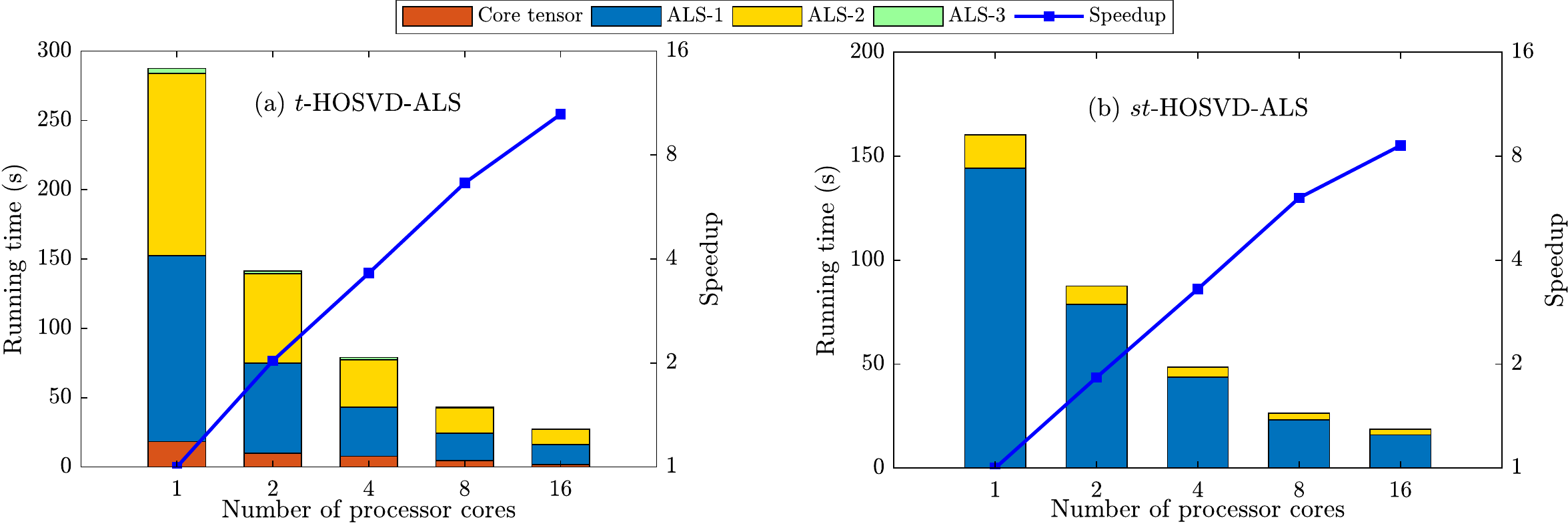}}
 \caption{\normalsize The running time and speedup of the ALS-based truncated HOSVD algorithms for compressing low multilinear rank tensor on a parallel computer.}\label{F6.4}
\end{figure}

An advantage of the proposed ALS-based methods is that they are easy to parallelize. To examine the parallel performance, in this experiment, 
we implement the $t$-HOSVD-ALS and $st$-HOSVD-ALS algorithms in C++ with OpenMP multi-threading parallelization \cite{OpenMP}.
The involved linear algebra operations are available with parallelization from the Intel MKL \cite{Intel,Wang2014} and the open-source ARMADILLO \cite{Sanderson2016,Sanderson2018} libraries. In the test, the input tensor is generated by \texttt{ttensor} in the Tensor Toolbox, whose size is $5000\times5000\times50$ and the multilinear rank is $(500,500,5)$. We break down the running time into different portions, including the ALS iterations along the three dimensions (ALS-$i$, $i\in\{1,2,3\}$) and the calculation of the core tensor. The test results are drawn in Fig. \ref{F6.4}. Also shown in the figures is the parallel scalability of the proposed algorithms.

From the figure, it can be seen that the proposed $t$-HOSVD-ALS and $st$-HOSVD-ALS can both scale well. 
In particular, $t$-HOSVD-ALS and $st$-HOSVD-ALS can achieve speedups of $10.50\times$ and $8.59\times$ as the number of processor cores is increased from 1 to 16, respectively.
Moreover, in both algorithms the ALS iterations, as the major costs, are accelerated efficiently with the increased number of processor cores. Another observation is that there is a slight drop of parallel efficiency when using 16 processor cores, which is caused by the fact that the factor matrix $\bm{U}^{(n)}$ is usually tall and skinny due to the low multilinear rank structure of tensor. Overall, the test results demonstrate that the proposed ALS-based algorithms are parallelization friendly and have good potential to scale further on larger high-performance computers.

\section{Conclusions}\label{S7}
In this paper, we proposed a class of ALS-based algorithms
for efficiently calculating the low multilinear rank approximation of tensors.
Compared with the original $t$-HOSVD and $st$-HOSVD algorithms,
the proposed algorithms are superior in several ways.
First, by eliminating the redundant computations of the singular vectors,
the overall costs of the algorithms are substantially reduced.
Second, the proposed algorithms are more flexible with adjustable convergence tolerance,
which is especially useful when the algorithms are used to generate initial solutions
for iterative methods such as HOOI.
Third, the proposed algorithms are free of the notorious date explosion issue
due to the fact that the ALS procedure does not explicitly require the intermediate matrices.
And fourth, the ALS-based approaches are parallelization friendly on high-performance computers.
Theoretical analysis shows that the ALS iteration in the proposed algorithms
is q-linear convergent with a relatively wide convergence region. 
Numerical experiments with both synthetic and real-world tensor data
demonstrate that proposed ALS-based algorithms can substantially reduce 
the total cost of low multilinear rank approximation and are highly parallelizable.

Possible future works could include applying of the proposed ALS-based algorithms to more applications,
among which we are especially interested in large-scale scientific computing.
It would also be of interest to study randomization techniques to further improve 
the performance of the proposed algorithms, considering the fact
that solving multiple least squares problems with different right-hand sides is the major cost. 
Some of the ideas presented in this work, such as the utilization of ALS for solving the intermediate 
low rank approximation problem, might be extended to other tensor decomposition models
such as tensor-train (TT) and hierarchical Tucker (HT) decompositions.

\section*{Acknowledgments} 
The authors would like to thank the anonymous reviewers for their valuable comments and suggestions that have greatly improved the quality of this paper.

\section*{Declarations}
\noindent\textbf{Funding} This study was funded in part by Guangdong Key R\&D Project (\#2019B121204008), Beijing Natural Science Foundation (\#JQ18001) and Beijing Academy of Artificial Intelligence.
\smallskip

\noindent\textbf{Conflicts of Interest} The authors declare that they have no conflict of interest.
\smallskip

\noindent\textbf{Availability of Data and Material} The datasets generated and analyzed during the current study are available from 
the corresponding author on reasonable request.
\smallskip

\noindent\textbf{Code Availability} The code used in the current study is available from the corresponding author on reasonable request.
\smallskip

\bibliographystyle{spmpsci}     
\bibliography{paper20-thosvd-als} 

\appendix
\section{Proof of Theorem \ref{thm2}} \label{A}

Based on the assumption of Theorem \ref{thm2},  $\bm{L}_{k}$ is nonsingular
and $\bm{L}_{k}^{T}\bm{L}_{k}$ is positive definite. Thus, the iterative form of Algorithm \ref{A3.3} is
\begin{equation}\label{E2.1}
\begin{split}
\bm{R}_{k} = \bm{A}^{T}\bm{L}_{k}(\bm{L}_{k}^{T}\bm{L}_{k})^{-1},\\
\bm{L}_{k+1} = \bm{AR}_{k}(\bm{R}_{k}^{T}\bm{R}_{k})^{-1},
\end{split}
\end{equation}
i.e.,
\begin{equation}\label{E2.2}
\bm{L}_{k+1} = \bm{AA}^{T}\bm{L}_{k}(\bm{L}_{k}^{T}\bm{AA}^{T}\bm{L}_{k})^{-1}(\bm{L}_{k}^{T}\bm{L}_{k}).
\end{equation}
Suppose that the full SVD of $\bm{A}$ is $\bm{A} = \bm{U}\bm{\Sigma} \bm{V}^{T}$, where
$$\bm{U} = [\bm{U}_{1},\bm{U}_{2}],\ \bm{V} = [\bm{V}_{1},\bm{V}_{2}], \,
\bm{\Sigma} = \left(
\begin{array}{cc}
\bm{\Sigma}_{1} & 0 \\
0 & \bm{\Sigma}_{2} \\
\end{array}
\right).
$$
Then from \eqref{E2.2}, we have
\begin{equation}\label{E2.3}
\bm{U}^{T}\bm{L}_{k+1} = \bm{U}^{T}\bm{AVV}^{T}\bm{A}^{T}\bm{UU}^{T}\bm{L}_{k}(\bm{L}_{k}^{T}\bm{UU}^{T}\bm{AVV}^{T}\bm{A}^{T}\bm{UU}^{T}\bm{L}_{k})^{-1}(\bm{L}_{k}^{T}\bm{UU}^{T}\bm{L}_{k}),
\end{equation}
which can be rewritten into block form
\begin{small}
	\begin{equation*}
	\left(
	\begin{array}{c}
	\bm{L}_{k+1}^{(1)} \\
	\bm{L}_{k+1}^{(2)} \\
	\end{array}
	\right) = \left(
	\begin{array}{cc}
	\bm{\Sigma}_{1}^{2} & 0 \\
	0 & \bm{\Sigma}_{2}\bm{\Sigma}_{2}^{T} \\
	\end{array}
	\right)\left(
	\begin{array}{c}
	\bm{L}_{k}^{(1)} \\
	\bm{L}_{k}^{(2)} \\
	\end{array}
	\right)
	(\bm{L}_{k}^{(1)T}\bm{\Sigma}_{1}^{2}\bm{L}_{k}^{(1)}+\bm{L}_{k}^{(2)T}\bm{\Sigma}_{2}\bm{\Sigma}_{2}^{T}\bm{L}_{k}^{(2)})^{-1}(\bm{L}_{k}^{(1)T}\bm{L}_{k}^{(1)}+\bm{L}_{k}^{(2)T}\bm{L}_{k}^{(2)}).
	\end{equation*}
\end{small}
It then follows that 
\begin{equation}\label{E2.5}
\begin{split}
\bm{L}_{k+1}^{(1)} &= \bm{\Sigma}_{1}^{2}\bm{L}_{k}^{(1)}(\bm{L}_{k}^{(1)T}\bm{\Sigma}_{1}^{2}\bm{L}_{k}^{(1)}+\bm{L}_{k}^{(2)T}\bm{\Sigma}_{2}\bm{\Sigma}_{2}^{T}\bm{L}_{k}^{(2)})^{-1}(\bm{L}_{k}^{(1)T}\bm{L}_{k}^{(1)}+\bm{L}_{k}^{(2)T}\bm{L}_{k}^{(2)}),\\
\bm{L}_{k+1}^{(2)} &= \bm{\Sigma}_{2}\bm{\Sigma}_{2}^{T}\bm{L}_{k}^{(2)}(\bm{L}_{k}^{(1)T}\bm{\Sigma}_{1}^{2}\bm{L}_{k}^{(1)}+\bm{L}_{k}^{(2)T}\bm{\Sigma}_{2}\bm{\Sigma}_{2}^{T}\bm{L}_{k}^{(2)})^{-1}(\bm{L}_{k}^{(1)T}\bm{L}_{k}^{(1)}+\bm{L}_{k}^{(2)T}\bm{L}_{k}^{(2)}).
\end{split}
\end{equation}
Furthermore, 
\begin{equation}\label{EA1}
\begin{split}
&(\bm{L}_{k}^{(1)T}\bm{\Sigma}_{1}^{2}\bm{L}_{k}^{(1)}+\bm{L}_{k}^{(2)T}\bm{\Sigma}_{2}\bm{\Sigma}_{2}^{T}\bm{L}_{k}^{(2)})^{-1} \\
=& (\bm{I}+(\bm{L}_{k}^{(1)T}\bm{\Sigma}_{1}^{2}\bm{L}_{k}^{(1)})^{-1}(\bm{L}_{k}^{(2)T}\bm{\Sigma}_{2}\bm{\Sigma}_{2}^{T}\bm{L}_{k}^{(2)}))^{-1}(\bm{L}_{k}^{(1)T}\bm{\Sigma}_{1}^{2}\bm{L}_{k}^{(1)})^{-1} \\
=& (\bm{I}+\sum\limits_{n=1}^{\infty}(-1)^{n}((\bm{L}_{k}^{(1)T}\bm{\Sigma}_{1}^{2}\bm{L}_{k}^{(1)})^{-1}(\bm{L}_{k}^{(2)T}\bm{\Sigma}_{2}\bm{\Sigma}_{2}^{T}\bm{L}_{k}^{(2)}))^{n})(\bm{L}_{k}^{(1)T}\bm{\Sigma}_{1}^{2}\bm{L}_{k}^{(1)})^{-1}.
\end{split}
\end{equation}
Here we suppose that the distance 
between $\mathcal{R}(\bm{L}_{k})$ and $\mathcal{R}(\bm{U}_{2})$ is small enough, 
therefore can be denoted as $\delta_{k}$, 
which only depends on $\|\bm{L}_{k}^{(2)}\|_{2}$ (i.e., there exist two constants $\alpha,\ \beta>0$ such that $\alpha\delta_{k}\leq\|\bm{L}_{k}^{(2)}\|_{2}\leq\beta\delta_{k}$). 
We can then obtain the lower bound of the distance between 
$\mathcal{R}(\bm{L}_{k})$ and $\mathcal{R}(\bm{U}_{1})$,
which is $\sqrt{1-\delta_{k}^{2}}$, and 
\begin{equation}\label{EA2}
\begin{split}
\|\bm{L}_{k}^{(1)}\|_{2}\leq C_{1}, \quad
\|(\bm{L}_{k}^{(1)})^{-1}\|_{2}\leq \frac{C_{2}}{\sqrt{1-\delta_{k}^{2}}},
\end{split}
\end{equation} 
where $C_{1},\ C_{2}$ are constants independent on $k$ and $\delta_{k}$. 
From \eqref{EA1} and \eqref{EA2}, there exists a constant $C$ that is only dependent on $C_{1},\ C_{2}$
so that the following inequality holds.
\begin{equation}\label{EA3}
(\bm{L}_{k}^{(1)T}\bm{\Sigma}_{1}^{2}\bm{L}_{k}^{(1)}+\bm{L}_{k}^{(2)T}\bm{\Sigma}_{2}\bm{\Sigma}_{2}^{T}\bm{L}_{k}^{(2)})^{-1}\leq(\bm{L}_{k}^{(1)T}\bm{\Sigma}_{1}^{2}\bm{L}_{k}^{(1)})^{-1}+\frac{C\delta_{k}^{2}}{(1-2\delta_{k}^{2})^{2}},
\end{equation}
Further, from \eqref{EA3} and \eqref{E2.5},  assume that  $\sigma_{r}>\sigma_{r+1}$,  we have
\begin{equation*}
\begin{split}
\bm{L}_{k+1}^{(2)} \leq \frac{\bm{\Sigma}_{2}\bm{\Sigma}_{2}^{T}}{\sigma_{r}^{2}}\bm{L}_{k}^{(2)}(\bm{L}_{k}^{(1)T}\frac{\bm{\Sigma}_{1}^{2}}{\sigma_{r}^{2}}\bm{L}_{k}^{(1)})^{-1}
(\bm{L}_{k}^{(1)T}\bm{L}_{k}^{(1)})
+\hat{C}\left(\frac{\delta_{k}^{2}}{(1-2\delta_{k}^{2})^{2}}+\frac{\delta_{k}^{2}}{1-\delta_{k}^{2}}+\frac{\delta_{k}^{4}}{(1-2\delta_{k}^{2})^{2}}\right),
\end{split}
\end{equation*}
where $\hat{C}$ is a constant. Clearly,
$$0<\bm{L}_{k}^{(1)T}\bm{L}_{k}^{(1)}\leq\bm{L}_{k}^{(1)T}\frac{\bm{\Sigma}_{1}^{2}}{\sigma_{r}^{2}}\bm{L}_{k}^{(1)},$$
and by Lemma \ref{lem1}, we have
$$\|(\bm{L}_{k}^{(1)T}\frac{\bm{\Sigma}_{1}^{2}}{\sigma_{r}^{2}}\bm{L}_{k}^{(1)})^{-1}(\bm{L}_{k}^{(1)T}\bm{L}_{k}^{(1)})\|_{2}\leq 1.$$
Since $\delta_{k}$ is small enough, we obtain
\begin{equation}\label{E2.7}
\begin{split}
\|\bm{L}_{k+1}^{(2)}\|_{2}\leq\frac{\sigma_{r+1}^{2}}{\sigma_{r}^{2}}\|\bm{L}_{k}^{(2)}\|_{2}+\tilde{C}\delta_{k}^{2}
\leq\frac{\sigma_{r+1}^{2}}{\sigma_{r}^{2}}\|\bm{L}_{k}^{(2)}\|_{2}+\frac{\tilde{C}}{\alpha^{2}}\|\bm{L}_{k}^{(2)}\|_{2}^{2},
\end{split}
\end{equation}
where $\alpha,\ \tilde{C}$ do not depend on $k$ and $\|\bm{L}_{k}^{(2)}\|_{2}$.

Denote 
$$
q = \frac{\sigma_{r+1}^{2}}{\sigma_{r}^{2}}+\frac{\tilde{C}}{\alpha^{2}}\|\bm{L}_{0}^{(2)}\|_{2}.
$$
Since we assume that $\mathcal{R}(\bm{L}_{0})$ is close to $\mathcal{R}(\bm{U}_{1})$ enough, $\|\bm{U}_{2}^{T}\bm{L}_{0}\|_{2}$ is sufficiently small, i.e., $\|\bm{L}_{0}^{(2)}\|_{2} = o(1)$. In other words, we assume that $q<1$. From (\ref{E2.7}), we have 
\begin{equation*}
\|\bm{L}_{k+1}^{(2)}\|_{2}\leq q\|\bm{L}_{k}^{(2)}\|_{2}
\end{equation*}
for all $k$, which leads to
$$\lim\limits_{k\rightarrow+\infty}\|\bm{L}_{k}^{(2)}\|_{2}\rightarrow0.$$
Combining with the assumption of $\bm{L}_{k}$, it is verified that $\mathcal{R}(\bm{L}_{k})$ is orthogonal to $\mathcal{R}(\bm{U}_{2})$ with $k\rightarrow+\infty$. Since the orthogonal complement space of $\mathcal{R}(\bm{U}_{2})$ is unique,  we have
$$\mathcal{R}(\bm{L}_{k}) = \mathcal{R}(\bm{U}_{1}),\ \ k\rightarrow+\infty,$$
where $\mathcal{R}(\bm{U}_{1})$ is the dominant subspace of $\bm{A}$. In other words, we have
\begin{equation*}
\begin{split}
\lim\limits_{k\rightarrow+\infty}\|\bm{L}_{k}\bm{L}_{k}^{\dag}-\bm{U}_{1}\bm{U}_{1}^{T}\|_{2}=0, 
\end{split}
\end{equation*}
where $\bm{L}_{k}^{\dag}$ is the pseudo-inverse of $\bm{L}_{k}$.
Further, from the iterative form of the ALS method, we have
$$\bm{R}_{k} = \bm{A}^{T}(\bm{L}_{k}^{\dag})^{T},$$
thus
$$\bm{L}_{k}\bm{R}_{k}^{T}=\bm{L}_{k}\bm{L}_{k}^{\dag}\bm{A}\rightarrow \bm{U}_{1}\bm{U}_{1}^{T}\bm{A},\ k\rightarrow+\infty.$$
And since the Frobenius norm $\|\cdot\|_{F}$ is continuous, 
$$\lim\limits_{k\rightarrow+\infty}\|\bm{A} - \bm{L}_{k}\bm{R}_{k}^{T}\|_{F}=\|\bm{A} - \bm{U}_{1}\bm{U}_{1}^{T}\bm{A}\|_{F}.$$
Since $\bm{U}_{1}\bm{U}_{1}^{T}\bm{A}$ is the exact solution of low rank approximation of $\bm{A}$,  the convergence of the ALS method is proved.

From \eqref{E2.7}, we further confirm the $q$-linear convergence of the ALS method, with approximate convergence ratio $\sigma_{r+1}^{2}/\sigma_{r}^{2}$.      \qed

\section{Proof of Theorem \ref{thm3}}\label{B}
The assumption of Theorem \ref{thm2} implies that $\bm{L}_{k}$ is nonsingular at every iteration $k$. We assume that $$\bm{L}_{k}^{(1)T}\bm{\Sigma}_{1}^{2}\bm{L}_{k}^{(1)}+\bm{L}_{k}^{(2)T}\bm{\Sigma}_{2}\bm{\Sigma}_{2}^{T}\bm{L}_{k}^{(2)}$$ 
is positive definite. Let $\varepsilon$ be a positive number such that $\sigma_{r}>\sigma_{r+1}-\varepsilon$. By (\ref{E2.5}), we know
\begin{equation}\label{E0}
\begin{split}
\bm{L}_{k+1}^{(1)} &= \bm{\Sigma}_{1}^{2}\bm{L}_{k}^{(1)}(\bm{L}_{k}^{(1)T}\bm{\Sigma}_{1}^{2}\bm{L}_{k}^{(1)}+\bm{L}_{k}^{(2)T}\bm{\Sigma}_{2}\bm{\Sigma}_{2}^{T}\bm{L}_{k}^{(2)})^{-1}(\bm{L}_{k}^{(1)T}\bm{L}_{k}^{(1)}+\bm{L}_{k}^{(2)T}\bm{L}_{k}^{(2)}),\\
\bm{L}_{k+1}^{(2)} &= \bm{\Sigma}_{2}\bm{\Sigma}_{2}^{T}\bm{L}_{k}^{(2)}(\bm{L}_{k}^{(1)T}\bm{\Sigma}_{1}^{2}\bm{L}_{k}^{(1)}+\bm{L}_{k}^{(2)T}\bm{\Sigma}_{2}\bm{\Sigma}_{2}^{T}\bm{L}_{k}^{(2)})^{-1}(\bm{L}_{k}^{(1)T}\bm{L}_{k}^{(1)}+\bm{L}_{k}^{(2)T}\bm{L}_{k}^{(2)}),
\end{split}
\end{equation}
which means
\begin{equation*}
\begin{split}
\bm{L}_{k+1}^{(1)} = \frac{\bm{\Sigma}_{1}^{2}}{(\sigma_{r}-\varepsilon)^{2}}\bm{L}_{k}^{(1)}(\bm{L}_{k}^{(1)T}\frac{\bm{\Sigma}_{1}^{2}}{(\sigma_{r}-\varepsilon)^{2}}\bm{L}_{k}^{(1)}+\bm{L}_{k}^{(2)T}\frac{\bm{\Sigma}_{2}\bm{\Sigma}_{2}^{T}}{(\sigma_{r}-\varepsilon)^{2}}\bm{L}_{k}^{(2)})^{-1}(\bm{L}_{k}^{(1)T}\bm{L}_{k}^{(1)}+\bm{L}_{k}^{(2)T}\bm{L}_{k}^{(2)}),\\
\bm{L}_{k+1}^{(2)} = \frac{\bm{\Sigma}_{2}\bm{\Sigma}_{2}^{T}}{(\sigma_{r}-\varepsilon)^{2}}\bm{L}_{k}^{(2)}(\bm{L}_{k}^{(1)T}\frac{\bm{\Sigma}_{1}^{2}}{(\sigma_{r}-\varepsilon)^{2}}\bm{L}_{k}^{(1)}+\bm{L}_{k}^{(2)T}\frac{\bm{\Sigma}_{2}\bm{\Sigma}_{2}^{T}}{(\sigma_{r}-\varepsilon)^{2}}\bm{L}_{k}^{(2)})^{-1}(\bm{L}_{k}^{(1)T}\bm{L}_{k}^{(1)}+\bm{L}_{k}^{(2)T}\bm{L}_{k}^{(2)}).
\end{split}
\end{equation*}
Clearly it holds that
\begin{equation}\label{Ethm2}
\begin{split}
\|\bm{L}_{k+1}^{(2)}\|_{2}\leq\frac{\sigma_{r+1}^{2}}{(\sigma_{r}-\varepsilon)^{2}}\|\bm{L}_{k}^{(2)}\|_{2}
\end{split}
\end{equation} 
under the condition that 
\begin{equation}\label{E0.0}
\bm{L}_{k}^{(1)T}\bm{L}_{k}^{(1)}+\bm{L}_{k}^{(2)T}\bm{L}_{k}^{(2)}\leq \bm{L}_{k}^{(1)T}\frac{\bm{\Sigma}_{1}^{2}}{(\sigma_{r}-\varepsilon)^{2}}\bm{L}_{k}^{(1)}+\bm{L}_{k}^{(2)T}\frac{\bm{\Sigma}_{2}\bm{\Sigma}_{2}^{T}}{(\sigma_{r}-\varepsilon)^{2}}\bm{L}_{k}^{(2)}.
\end{equation}
If $\bm{L}_{k}^{(1)}$ is nonsingular, then \eqref{E0.0} implies 
\begin{equation}\label{E1}
(\bm{L}_{k}^{(2)}(\bm{L}_{k}^{(1)})^{-1})^{T}(\bm{I}-\frac{\bm{\Sigma}_{2}\bm{\Sigma}_{2}^{T}}{(\sigma_{r}-\varepsilon)^{2}})(\bm{L}_{k}^{(2)}\bm{L}_{k}^{(1)-1})\leq \frac{\bm{\Sigma}_{1}^{2}}{(\sigma_{r}-\varepsilon)^{2}}-\bm{I}.
\end{equation}
It follows to see that 
\begin{equation}\label{E2}
\|\bm{L}_{k}^{(2)}(\bm{L}_{k}^{(1)})^{-1}\|_{2}\leq\sqrt{\frac{\sigma_{r}^{2}-(\sigma_{r}-\varepsilon)^{2}}{(\sigma_{r}-\varepsilon)^{2}-\sigma_{min}^{2}}}
\end{equation}
is a sufficient condition of \eqref{E1}.

Next we will prove that if the initial guess $\bm{L}_{0}$ satisfies condition \eqref{E2}, then $\bm{L}_{k}^{(1)}$ is nonsingular and \eqref{Ethm2} is satisfied at every iteration $k$. 

Provided that $\bm{L}_{0}$ satisfies \eqref{E2},  we obtain
\begin{equation}\label{E3}
\|\bm{L}_{1}^{(2)}\|_{2}\leq\frac{\sigma_{r+1}^{2}}{(\sigma_{r}-\varepsilon)^{2}}\|\bm{L}_{0}^{(2)}\|_{2}.
\end{equation}
And according to the proof of Theorem \ref{thm2}, we know $\bm{L}_{1}^{(1)}$ is also nonsigular, which implies $$\bm{L}_{1}^{(1)T}\bm{\Sigma}_{1}^{2}\bm{L}_{1}^{(1)}+\bm{L}_{1}^{(2)T}\bm{\Sigma}_{2}\bm{\Sigma}_{2}^{T}\bm{L}_{1}^{(2)}$$ is positive definite. Then by \eqref{E0}, we have  
\begin{equation}\label{E5}
\|\bm{L}_{1}^{(2)}(\bm{L}_{1}^{(1)})^{-1}\|_{2} \leq \|\bm{\Sigma}_{2}\bm{\Sigma}_{2}^{T}\bm{L}_{0}^{(2)}(\bm{L}_{0}^{(1)})^{-1}\bm{\Sigma}_{1}^{-2}\|_{2}\leq\frac{\sigma_{r+1}^{2}}{\sigma_{r}^{2}}\|\bm{L}_{0}^{(2)}(\bm{L}_{0}^{(1)})^{-1}\|_{2}\leq\sqrt{\frac{\sigma_{r}^{2}-(\sigma_{r}-\varepsilon)^{2}}{(\sigma_{r}-\varepsilon)^{2}-\sigma_{min}^{2}}}.
\end{equation}
Analogously, we can prove that for every iteration $k$, $\bm{L}_{k}^{(1)}$ is nonsingular, i.e., $\bm{L}_{k}^{(1)T}\bm{\Sigma}_{1}^{2}\bm{L}_{k}^{(1)}+\bm{L}_{k}^{(2)T}\bm{\Sigma}_{2}\bm{\Sigma}_{2}^{T}\bm{L}_{k}^{(2)}$ is positive definite, and \eqref{E2} is satisfied. Since \eqref{E2} is a sufficient condition of \eqref{E1} and \eqref{E0.0}, inequality \eqref{Ethm2} is true at every iteration $k$, which implies that
$$\lim\limits_{k\rightarrow0}\|\bm{L}_{k}\bm{L}_{k}^{\dagger}-\bm{U}_{1}\bm{U}_{1}^{T}\|_{2}=0.$$  
The rest part of the proof is analogous to the proof of Theorem \ref{thm2}, which is omitted for brevity.  \qed

\section{Proof of Theorem \ref{thm4}} \label{C}

In Algorithm \ref{A3.4}, $\bm{U}^{(n)}$ is obtained from the rank-$R_{n}$ approximation of $\bm{A}_{(n)}$,
which is done in an iterative manner and allows a tolerance parameter $\eta_{n}$. Therefore, we have  
\begin{equation}\label{E5.3}
\|\bm{A}_{(n)} - \bm{L}^{*}\bm{R}^{*T}\|_{F}^{2}\leq \eta_{n}^{2}\|\bm{\mathcal{A}}\|_{F}^{2} + \gamma_{n},
\end{equation}
where $\bm{L}^{*}$ and $\bm{R}^{*}$ are the same as in Algorithm \ref{A3.4}. 
Note that $\bm{R}$ is updated by solving a multi-side least squares problem
\begin{equation*}
\min\limits_{\bm{R}}\|\bm{A}_{(n)} - \bm{LR}^{T}\|_{F},
\end{equation*}
whose exact solution is $\bm{R} = \bm{A}_{(n)}^{T}(\bm{L}^{\dagger})^{T}$. Thus 
\begin{equation}\label{E5.5}
\|\bm{A}_{(n)} - \bm{L}^{*}\bm{R}^{*T}\|_{F} = \|\bm{A}_{(n)} - \bm{L}^{*}\bm{L}^{*\dagger}\bm{A}_{(n)}\|_{F},
\end{equation}
where $\bm{L}^{*}\bm{L}^{*\dagger}$ represents an orthogonal projection on subspace $\mathcal{R}(\bm{L}^{*})$. 
Consequently, by \eqref{E5.3}, \eqref{E5.5} and \eqref{E5.1}, we have
\begin{equation*}
\|\hat{\bm{\mathcal{A}}} - \bm{\mathcal{A}}\|_{F}^{2}\leq\sum\limits_{n=1}^{N}(\eta_{n}^{2}\|\bm{A}_{(n)}\|_{F}^{2}+\gamma_{n}),
\end{equation*}
which means
\begin{equation*}
\frac{\|\hat{\bm{\mathcal{A}}} - \bm{\mathcal{A}}\|_{F}^{2}}{\|\bm{\mathcal{A}}\|_{F}^{2}}\leq\sum\limits_{n=1}^{N}(\eta_{n}^{2}+\frac{\gamma_{n}}{\|\bm{\mathcal{A}}\|_{F}^{2}}).
\end{equation*}
Combining with \eqref{E5.1}, we obtain \eqref{E5.62}.

The error analysis of Algorithm \ref{A3.5} can be analogously done.  \qed

\end{document}